\documentclass[12pt]{article}
\usepackage{latexsym,amssymb,a4,epsfig}

\newcommand{\ES}{\mbox{\rlap{$\bigcirc$}\kern2.5pt/\kern2.5pt}}
\newcommand{\NS}{\mbox{\rlap{I}\kern1.5pt\rlap{I}\kern0.4ptN}}
\newcommand{\ZS}{\mbox{\rlap{Z}\kern1.5ptZ}}
\newcommand{\RS}{\mbox{\rlap{I}\kern1.5ptR}}
\newcommand{\QS}{\mbox{\rlap{Q}\kern3pt\vrule height 6.6pt depth 0pt\kern4.777pt}}
\newcommand{\CS}{\mbox{\rlap{C}\kern3pt\vrule height 6.6pt depth 0pt\kern4.222pt}}
\newcommand{\curlR}{{\cal R}}

\newcommand{\curlC}{{\cal C}}
\newcommand{\curl}[1]{{\cal #1}}

\newcommand{\defn}[1]{{\it #1}}
\newcommand{\qed}{\hspace*{\fill}$\epsfxsize.08in\epsffile{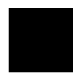}$\\}

\newenvironment{proof}
 {\begin{trivlist} \item[\hskip \labelsep {\small\bf Proof.\phantom{three}}]
		   \nopagebreak}
 {\nopagebreak
    \hfill$\raisebox{-2mm}{\epsfxsize.08in\epsffile{proof.ps}}$\end{trivlist}}

\newenvironment{sketchproof}
 {\begin{trivlist} \item[\hskip \labelsep {\small\bf Sketch proof.\phantom{three}}]
		    \nopagebreak}
 {\nopagebreak
    \hfill$\raisebox{-2mm}{\epsfxsize.08in\epsffile{proof.ps}}$\end{trivlist}}

\newcommand{\leftpic}{\epsffile{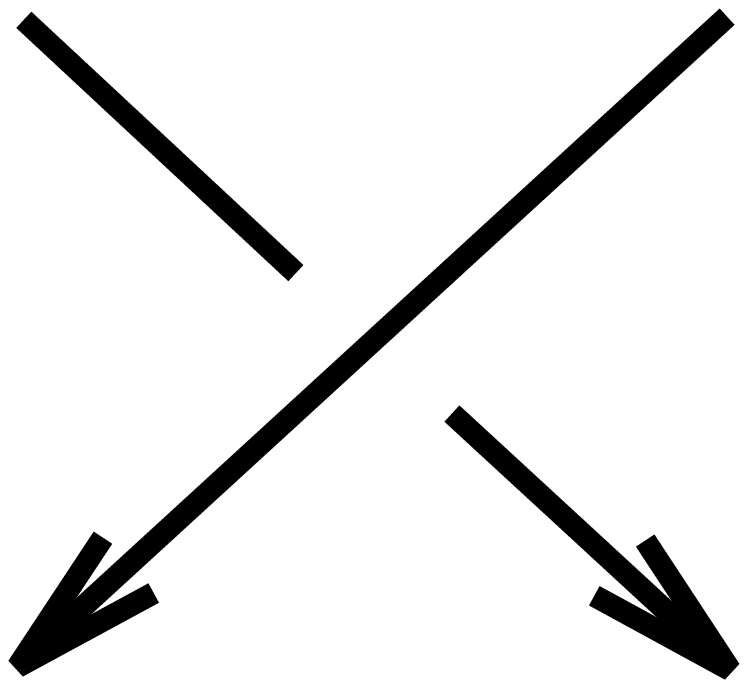}}
\newcommand{\rightpic}{\epsffile{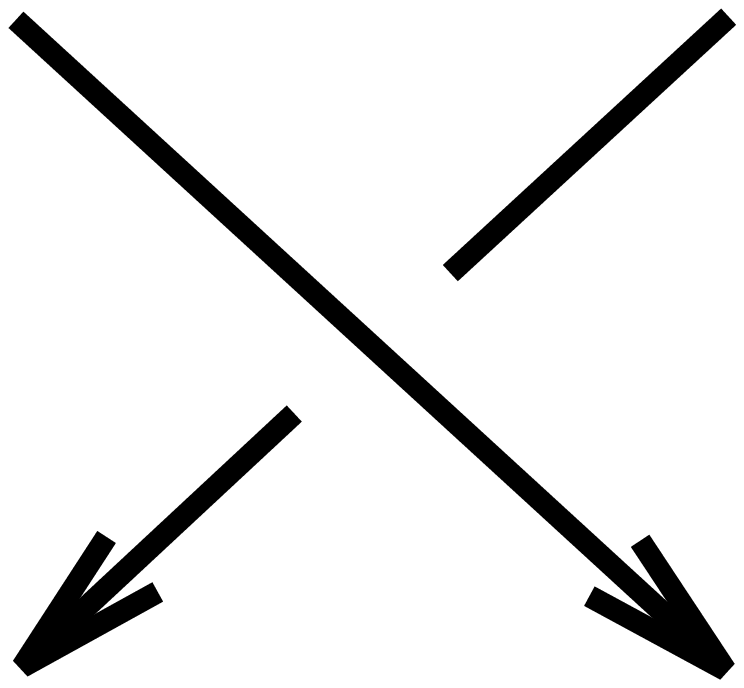}}
\newcommand{\framel}{\epsffile{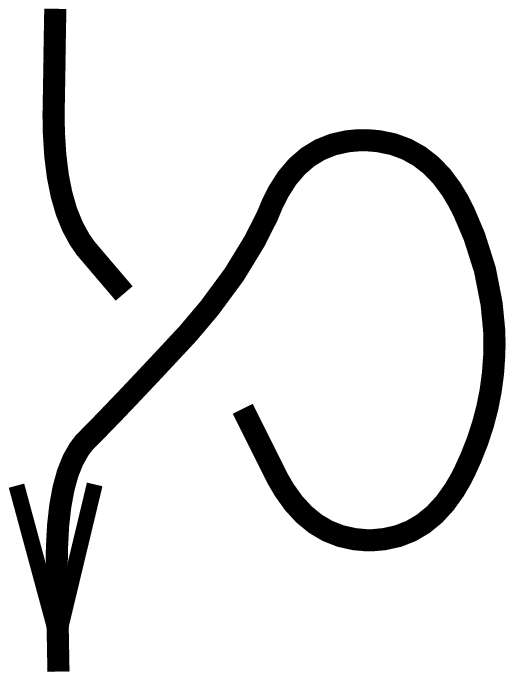}}
\newcommand{\orline}{\epsffile{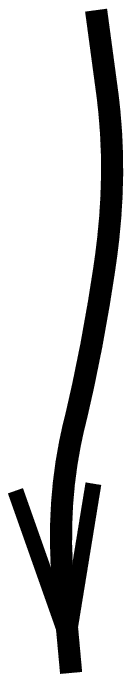}}
\newcommand{\parrallel}{\epsffile{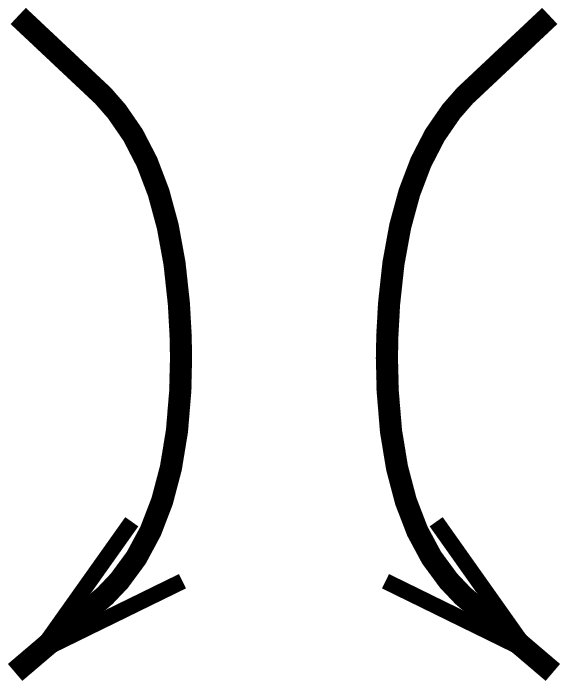}}

\begin{document}
\bibliographystyle{bibsty}

\title{Adams operators and knot decorations.}
\author{ 
A.K. Aiston\thanks{Supported by a EPSRC grant GR/J72332.}\\
Dept. Math. Sci., Liverpool University, Liverpool, L69 3BX.}
\maketitle
\pagestyle{myheadings}
\markright{Adams Operators and Knot Decorations.}

\begin{abstract}

We use an explicit isomorphism from the representation ring of the quantum
group $U_q(sl(N))$
to the Homfly skein of the annulus, to determine an element of the skein
which is the image of $\psi_m(c_1)$,  the $m$th Adams operator on the 
fundamental representation. This element is a linear combination of $m$ very
simple $m$-string braids.

Using this skein element, we show that the Vassiliev invariant of degree
$n$ in the power series expansion of the $U_q(sl(N))$ 
quantum invariant of a knot coloured
by $\psi_m(c_1)$ is the canonical Vassiliev invariant with weight system
$W_n\circ \psi_m^{(n)}$, where $W_n$ is the weight system for the Vassiliev
invariant of degree $n$ in the expansion of the quantum invariant of the knot
coloured by $c_1$ and $\psi_m^{(n)}$ is the Adams operator
on $n$-chord diagrams defined by 
Bar-Natan \cite{barnatan}.
\end{abstract}

\subsection{Introduction.}

A motivation for this article and \cite{qdim,idemp} was to investigate
the extent to which the properties of the Homfly polynomial and the 
$U_q(sl(N))$-quantum invariants are determined by the combinatoric 
properties of Young diagrams.
This extends the work of Morton \cite{nato}.
The central role of Young diagrams in both the study of 
symmetric polynomials and the representation theory of Lie algebras 
is described in many texts, for example \cite{macnew,repthry}.

Quantum groups, introduced by Drinfel'd \cite{drin2}, are $1$-parameter
deformations of classical Lie algebras.  We are concerned in particular
with the quantum group $U_q(sl(N))$.   
For generic $q$, Lusztig \cite{lusz} showed that the
representation theory of $U_q(sl(N))$ is isomorphic to that of 
the classical Lie algebra. Thus the Young diagrams play a role in
the world of quantum groups.

The Homfly polynomial \cite{homfly,pt} is an oriented link invariant which
can be described combinatorially in terms of skein relations.
We will work with a framed $3$-variable version, ${\cal X}(x,v,s)$.  
At certain evaluations of the variables, Turaev \cite{tur}
showed $\cal X$ is the $U_q(sl(N))$-invariant of a knot
coloured by the fundamental representation. 
Further, there are patterns
in the Homfly skein of the annulus for which the Homfly 
polynomial of the satellite of a knot, at these values, is the quantum 
invariant of the knot coloured by some, higher dimensional, representation.
Since Young diagrams index these representations, they play an integral role 
in Homfly skein theory.

The layout of this paper is as follows.
Section \ref{sec2} recalls some facts about Young diagrams.
Section \ref{sec3} discusses Homfly skein theory.
In Sect.~\ref{sec4} we recall the relevant results of \cite{idemp,mine}.
Section \ref{sec5} constructs the isomorphism between the representation ring
of the quantum groups and the Homfly skein of the annulus.
In Sect.~\ref{sec6} the $m$th Adams operator of the fundamental representation,
$\psi_m(c_1)$ is presented.  Our main result, Theorem~\ref{xbiff},
gives the image of $\psi_m(c_1)$ under the isomorphism of Sect.~\ref{sec5}
as a linear combination of $m$ very simple $m$-string braids.
In Sect.~\ref{sec7} we prove that the Adams operators give rise to 
canonical Vassiliev invariants whose weight systems are defined using
the chord diagram Adams operators of Bar-Natan\cite{barnatan}.

Throughout, $\Lambda=\QS[x^{\pm 1}, v^{\pm 1}, s^{\pm 1}]$
will be the ring of Laurent polynomials in $x$, $v$ and $s$ 
and we will denote $(s^i-s^{-i})/(s-s^{-1})$ by $[i]$.

\subsection{The ring of Young diagrams.}
\label{sec2}

There is a wealth of detail about the features of Young
diagrams in many texts (for example \cite{repthry, macnew}).
Many of the properties to the fore in this paper are discussed in 
\cite{idemp} and we will not repeat them here.
However, we use an example to fix notation.

\begin{figure}[ht]
\[
T(\nu)\ =\ \raisebox{-3mm}{\epsfxsize.533in\epsffile{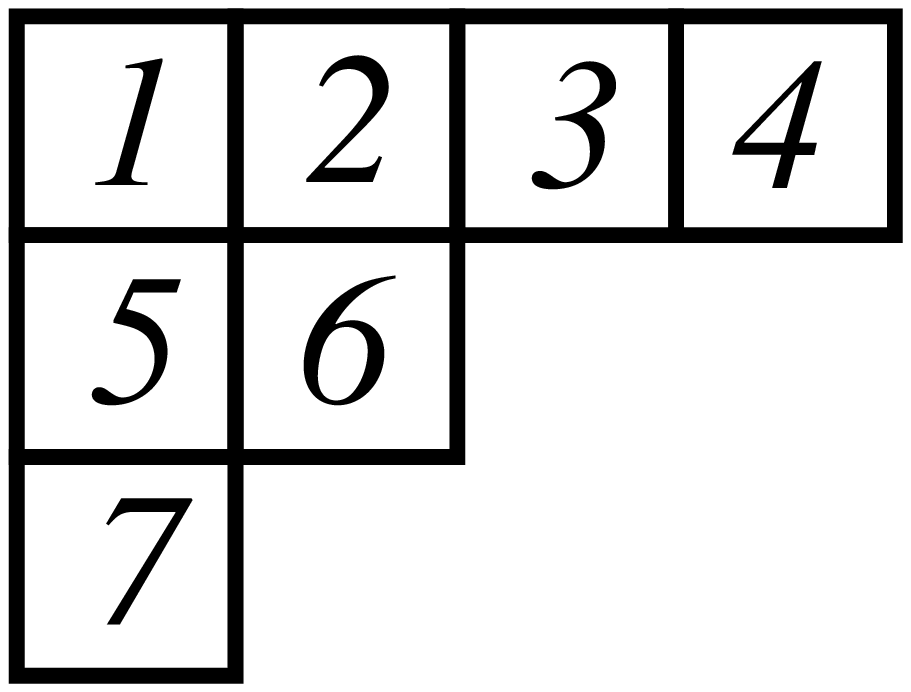}}
\]
\caption{The tableau $T(\nu)$.}
\label{eaux}
\end{figure}

Let $\nu$ be the Young diagram $(4,2,1)$. The transpose of $\nu$ 
is $\nu^\vee=(3,2,1,1)$. The Young tableau $T(\nu)$
is the assignment of numbers to the cells of $\nu$
shown Fig.~\ref{eaux} and the permutation $\pi_\nu$
is defined as $\pi_\nu(i)=j$
where transposition takes the cell $i$ in $T(\nu)$ 
to the cell $j$ in  $T(\nu^\vee)$.
Thus $\pi_\nu=\left(2\,4\,7\,3\,6\,5\right)$. 

The set of formal 
$\Lambda$-linear combinations of Young diagrams can be given an
associative, commutative algebra structure.  
The identity is the partition $(0)$
and the structure constants are taken to be the 
Littlewood-Richardson  coefficients $\{a_{\lambda\mu}^\nu\}$.
These can be defined as follows.
Let $\lambda=(\lambda_1,\cdots,\lambda_k)$
and $\mu=(\mu_1,,\cdots,\mu_m)$ be two Young diagrams. 
A $\mu$-\defn{expansion} of $\lambda$ is obtained by first  
adding $\mu_1$ boxes to $\lambda$, each labelled with a $1$.
No two boxes can be placed in the same column and the result must be
a legitimate Young diagram.  Then add $\mu_2$ boxes
labelled $2$ (respecting the same rules) and continue until
you have added $\mu_m$ boxes labelled $m$.  
At each stage no two cells with the same label can appear in
the same column.
For any given cell, let $n_i$ be the number of cells numbered 
$i$ above and to
the right of it (including the cell itself).
\noindent The expansion is called \defn{strict} if, 
for every cell, $n_i\geq n_j$ whenever $i<j$. 
The coefficient $a_{\lambda\mu}^\nu$ is the number of ways $\nu$
can be obtained from $\lambda$ by a strict $\mu$-expansion.
This coefficient is also the multiplicity of the 
irreducible $U_q(sl(N))$-representation indexed by $\nu$ in the decomposition 
of the tensor product of the representations indexed by $\lambda$ and 
$\mu$, however, our definition is independent of $N$.

It is well known (see for example \cite{macnew}) that the 
algebra of Young diagrams $Y$ is freely generated as a 
polynomial algebra by the Young diagrams with a single column.
If we denote the $\Lambda$-algebra 
of polynomials in an infinite 
number of indeterminates, $\{c_i\}_{i \in {\mathbb N}}$, 
by $\curlR_\infty$
then, since they are both freely generated on
a countably infinite set of generators, $\curlR_\infty$ is
isomorphic to $Y$, via the algebra isomorphism
\[
\phi\ :\ c_i\quad\longmapsto\quad\raisebox{-5mm}{\epsfysize.5in\epsffile{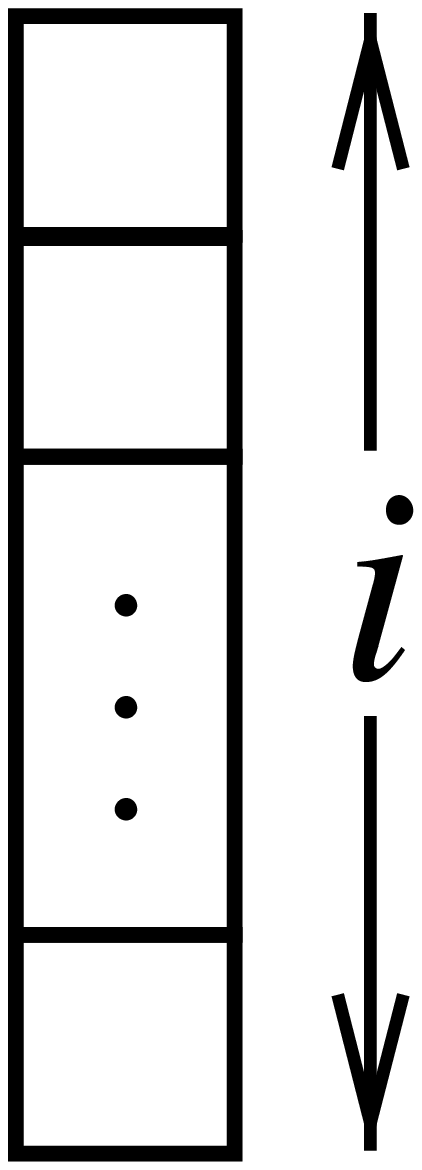}}\,.
\]
\label{move}

Let ${\cal R}_N$ denote the representation ring of the classical complex 
Lie algebra $sl(N)$ (and thus by Lusztig also the representation ring of 
the quantum group $U_q(sl(N))$ for generic $q$). 
The elements of the ring are indexed by the elements of $Y$
and the following relations hold: 
\newcounter{anna}
\begin{list}{\roman{anna})}{\usecounter{anna}}
        \setlength{\rightmargin}{\leftmargin}  
\item A representation indexed by a Young diagram with a column containing 
more than $N$ cells is $0\in{\cal R}_N$.
\item A representation indexed by a Young diagram with a column containing
exactly $N$ cells is isomorphic to the representation indexed by the diagram
with that column removed.
\end{list}  
Thus,  we have an algebra homomorphism, 
$p_N:{\cal R}_\infty \rightarrow {\cal R}_N$,
taking $c_k$ to the $k$th exterior power of
the fundamental representation.  
Any element of ${\cal R}_\infty$, which is sent, by $\phi$ to a single Young 
diagram, is mapped under $p_N$ to an irreducible 
representation. The set of Young diagrams with 
fewer than $N$ rows is in one-to--one correspondence with the irreducible 
representations.  Thus, ${\curlR}_N$ is isomorphic to 
${\cal R}_\infty/\left<c_k \;\forall k>N\;\quad c_N-1\right>$
and generated by $\{c_1\cdots c_{N-1}\}$.

Let $d_l$ denote the element of ${\cal R}_\infty$ whose image under $\phi$
is $(l)$.  Then $p_N(d_l)$ is the $l$th symmetric power of $c_1$.
In particular $d_1=c_1$ and $d_0=c_0=1$.

\begin{figure}
\[
\raisebox{-5mm}{\epsfxsize0.5in\epsffile{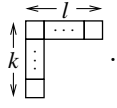}}\,.
\]
\caption{The Young diagram $\mu_{k,l}$.}
\label{hook}
\end{figure}

Consider the strict expansions of 
a Young diagram with a single column by one with a single row.
If we add the cells from the single row to 
the single column each cell must be placed in
a different column.  Hence, there are only
two possibilities;  one cell is added to the first column
and all the others start new columns or all the cells from
the single row start a new column.  The corresponding
diagrams have the shape of a hook.  Since $\phi$
is an isomorphism, we shall use $c_k$ and $d_l$ to
denote both the elements of $\curl{R}_\infty$ and their associated
Young diagrams.  For $k,l>0$, write $\mu_{k,l}\in Y$ for the Young 
diagram in Fig.~\ref{hook}, with $k+l-1$ cells.
Then,
\begin{equation}
c_kd_l=\left\{ \begin{array}{cl} 
                     c_k& \mbox{ if $l=0$}\\
                     d_l& \mbox{ if $k=0$}\\
                     \mu_{k+1,l}+\mu_{k,l+1}& \mbox{ otherwise.}
                   \end{array}
           \right.
\label{ckdl}
\end{equation}
Note that \,$\mu_{k,1}=c_k$ and that \,$\mu_{1,l}=d_l$.
The following classical result plays
an integral role in what follows.

\subsubsection{Proposition.{\rm\protect\cite{weyl}}}
\label{equal1}

Let $C(X)=\displaystyle\sum_{k=0}^{\infty}(-1)^k c_kX^k$ and 
$D(X)=\displaystyle\sum_{l=0}^{\infty}d_l X^l$ be formal power series 
with coefficients in $Y$.  These series satisfy the relation
$C(X)D(X)=1$.
\begin{proof}
Let $C(X)D(X)=\sum_{m=0}^{\infty}a_mX^m$ with
$a_m=\sum^m_{k=0} (-1)^kc_kd_{m-k}\,.$
Then $a_0=c_0d_0=1$.
For $m>0$,
\begin{eqnarray*}
      a_m&=&\sum_{k=0}^m(-1)^kc_kd_{m-k}\\
         &=&d_m+\sum_{k=1}^{m-1}\left((-1)^k(\mu_{k+1,m-k}+\mu_{k,m-k+1})\right)
				\quad+\quad(-1)^mc_m \\
         &=&d_m+\sum_{k=2}^{m}(-1)^{k-1}\mu_{k,m-k+1}
                \quad+\quad\sum_{k=1}^{m-1}(-1)^k\mu_{k,m-k+1}\quad+\quad(-1)^mc_m\\
         &=&d_m+(-1)^{m-1}\mu_{m,1}+(-1)\mu_{1,m}+(-1)^mc_m\\
         &=&d_m-d_m+(-1)^{m-1}(c_m-c_m)\\
         &=&0.
\end{eqnarray*}
\end{proof}
Let $C_N(X)$ denote the polynomial obtained from $C(X)$ by taking 
coefficients to lie in $\curl{R}_N$.  Thus,
$C_N(X)=1-c_1X+\cdots+(-1)^N c_N X^N$.
This polynomial can be formally factorised as 
$C_N(X)=\displaystyle\prod_{i=1}^N(1-x_iX)$ and each 
$c_k\in\curl{R}_N$ can be written, formally, as the $k$th elementary symmetric 
function in  $\{x_i\}_{i=1}^N$.
For example, $c_1=\displaystyle\sum_{i=1}^N x_i$.
The algebra ${\cal R}_N$ can, therefore, be presented as the 
quotient of the algebra of
symmetric polynomials in $N$ variables by the ideal generated by the 
$c_N-1$.  However, for the purposes of this paper
if we wish to work with a particular degree of polynomial, 
we can choose $N$ large enough so that we do not need to be concerned with
this extra relation.
Under this interpretation, $d_l$ becomes the $l$th complete symmetric 
polynomial
(the sum of all the monomials of degree $l$).
Both the set of elementary symmetric polynomials and the complete symmetric 
polynomials generate the algebra of symmetric polynomials.  In Sect.~\ref{sec4}
we shall give skein theoretic version of both these generating sets.  We then 
go on to give skein theoretic version of a third generating set, 
the power sums.

\subsection{The Homfly skein.}
\label{sec3}

We give a brief description of skein theory based on planar pieces
of knot-diagrams and a framed version of the Homfly polynomial.
The ideas go back to Conway and have been substantially developed
by Lickorish and others.  A fuller version of this account can be
found in \cite{nato, idemp}. 

We shall work with the framed Homfly polynomial
${\cal X}$. This is an invariant of framed oriented links,
 constructed from the Homfly polynomial by setting
${\cal X}(L)=(xv^{-1})^{\omega (D)}P(L)$, where $\omega (D)$
is the writhe of any diagram $D$ of the framed link $L$ which 
realises the chosen framing by means of the `blackboard parallel'. 
With the normalisation that ${\cal X}$ takes the value $1$ on the
empty knot, ${\cal X}$ is uniquely determined by the skein relations
in Fig.~\ref{quadratic}.
\begin{figure}[ht]
\[
x^{-1}\,{\cal X}\left(\raisebox{-1mm}{\,\epsfxsize.2in\leftpic\,}\right)
        \,-\, x\,{\cal X}\left(\raisebox{-1mm}{\,\epsfxsize.2in\rightpic\,}
                                                                        \right)
        \,=\,z\,{\cal X}\left(\raisebox{-2mm}{\,\epsfxsize0.2in\parrallel\,}
                                                        \right) \,,
\qquad
{\cal X}\left(\raisebox{-2mm}{\,\,\epsfxsize0.2in\framel\,}\right)
        \ =\ (xv^{-1})\,{\cal X}\left(
                        \raisebox{-3mm}{\,\,\epsfxsize0.04in\orline\,\,}
                                        \right)\,.
\]
\caption{The framed Homfly skein relations.}
\label{quadratic}
\end{figure}

Let $F$ be a planar surface and consider diagrams in $F$ consisting of oriented arcs joining 
any distinguished boundary points and oriented closed curves, 
up to regular isotopy. They carry the
implicit framing defined by the parallel curves in the diagram. Define the
{\it framed Homf\,ly skein\/} of $F$, denoted by ${\cal S}(F)$, to consist of
linear combinations of diagrams in $F$ modulo
the skein relations in Fig.~\ref{quadratic}.
A detailed description of the general theory appears in \cite{nato}
but here we are interested in  three specific cases. 

When $F$ is the whole  plane
${\mathbb R}^2$ then ${\cal S}({\mathbb R}^2$) is just the set of
linear combinations of
framed link diagrams, modulo the skein relations. Every diagram $D$ represents
a scalar multiple, namely  ${\cal X}(D)$, of the empty diagram.

Let $R^n_n \cong I\times I$ be the rectangle with $n$ distinguished  points
on its top and bottom edge.
We insist
that any arcs in $R^n_n$ enter at the top and leave at the bottom.
Diagrams in $R^n_n$ are termed {\it oriented\/} $n$-{\it tangles\/}, and 
include the case of
$n$-string braids.

Originally defined by Elrifai and Morton \cite{fai},
for each permutation $\pi\in S_n$ there is an $n$-string
{\it positive permutation braid\/} (which we shall abbreviate to p.p.b.), 
$\omega_\pi$,
uniquely determined by the fact that it is the minimal length braid, 
with all crossings positive, which joins the point $i$ at the top to $\pi(i)$ 
at the bottom. 
The {\it negative permutation braid\/}(n.p.b.) $\overline{\omega}_\pi$ is
defined in exactly the same manner, except that all the crossing are negative. 
Morton and Traczyk \cite{mt} proved that the 
the $n!$ $n$-string p.p.b.s 
are a linear basis for ${\cal S}(R_n^n)$.

The skein forms an algebra over $\Lambda$ with
multiplication derived from the concatenation of diagrams: write $ST$ 
for the diagram obtained by placing  $S$ above $T$.
The resulting algebra is a quotient of the braid-group algebra, shown
in \cite{mt} to be isomorphic to the Hecke algebra $H_n$ of type $A$, 
with the explicit presentation
\[
H_n\quad=\quad
\left<
\begin{array}{ccc}
        \begin{array}{ccc}
        \sigma_i& : & i=1,\ldots, n-1\\
               &   &
        \end{array}
&\left.\begin{array}{c}
         \\
         \\ 
        \end{array}     \right\vert

&       \begin{array}{l}
        \sigma_i\sigma_j=\sigma_j\sigma_i~:~\vert i-j\vert>1\\
        \sigma_i\sigma_{i+1}\sigma_i=\sigma_{i+1}\sigma_i\sigma_{i+1}\\
        x^{-1}\sigma_i-x\sigma^{-1}_i=z\;,
        \end{array}
\end{array}
\right>\;,
\]
where $\sigma_i$ denotes the p.p.b. with permutation
$(i\;i+1)$.

A {\it wiring\/} $W$ of a surface $F$ into another surface $F'$
is a choice of inclusion of $F$ into
$F'$ and a choice of a fixed diagram of curves and arcs in $F'- F$ whose
 boundary  is the union of the distinguished
sets of $F$ and $ F'$. A wiring $W$ determines naturally a $\Lambda$-linear map
${\cal S}(W):{\cal S}(F)\to {\cal S}(F')$.

\begin{figure}[ht]
\[
\epsfysize.6in\epsffile{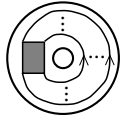}
\]
\caption[]{The wiring of ${\cal S}(R_n^n)$ into ${\cal S}(S^1\times I)$}
\label{wire}
\end{figure}

We can wire the rectangle $R^n_n$
into the annulus as indicated in Fig.~\ref{wire}, this being our third surface.
The resulting diagram in the annulus is called the {\it closure\/} 
of the oriented
tangle.  We shall also use the term `closure' for the family of
$\Lambda$-linear maps from the Hecke algebras $H_n$ to the skein of 
the annulus induced by this wiring.

The skein of the annulus, ${\cal S}(S^1\times I)$, itself forms an algebra, the
product given by stacking the annuli one inside the other.  This is
obviously commutative (lift the inner annulus up and stretch it so
that the outer one will fit inside it).
Write ${\cal C}$ for ${\cal S}(S^1\times I)$ regarded as a  $\Lambda$-algebra
in this way.
Turaev \cite{turbas} showed that ${\cal C}$ is freely generated
as an algebra by $\{ A_m$, $m\in {\mathbb Z}\}$,  
where $ A_m$ is the closure of the p.p.b. for 
the cycle $(1\,2\,\ldots\,\vert m\vert)$ 
and depending on whether $m>0$ or $m<0$, 
we take the closure to be oriented in the same or the opposite sense 
to the core of the annulus.
(For an alternative account see Hoste and Kidwell \cite{hk}.)
For example
\[
A_2\quad=\quad\raisebox{-5mm}{\epsfxsize.6in\epsffile{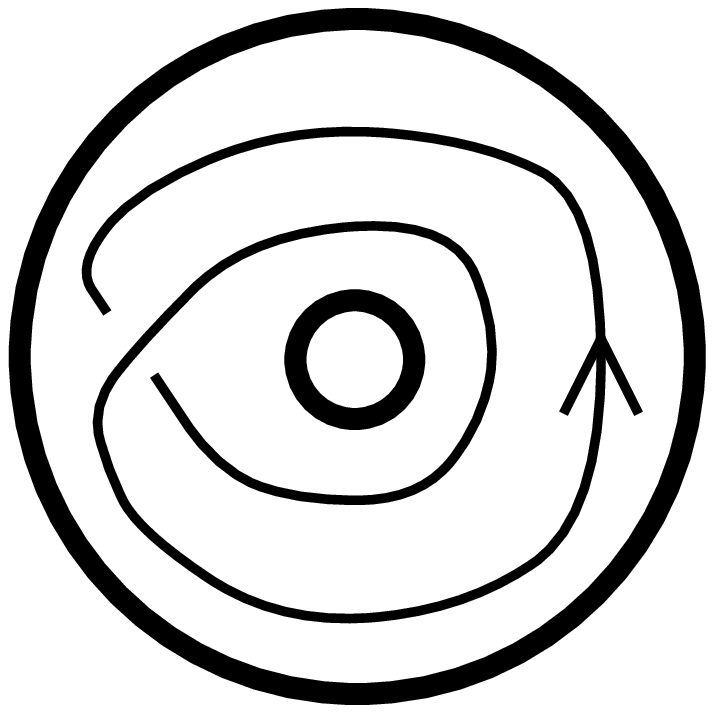}}
\qquad\mbox{ and }\qquad
A_{-2}\quad=\quad\raisebox{-5mm}{\epsfxsize.6in\epsffile{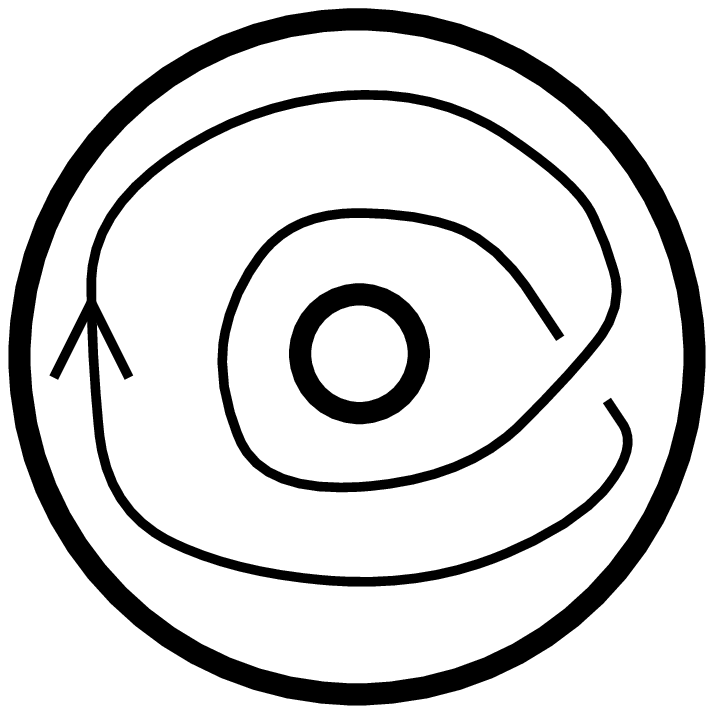}}\ \,.
\]
The identity element of the algebra, $A_0$, is 
represented by the empty diagram. 
\begin{figure}[ht]
\[
\epsfxsize.5in\epsffile{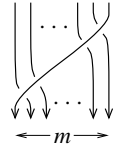}
\]
\caption{The braid which closes to $A_{\pm m}\in\curl{C}$.}
\label{am}
\end{figure}
The sub-algebra ${\cal C}^+$, spanned by the closures of
oriented tangles using the wiring shown in Fig.~\ref{wire},
\label{bazza}
is freely generated by $\{A_m\,:\,m \geq 0\}$.  
Set $\curl{C}^{(n)}$ to be the linear space generated by
all closures of $n$-string oriented tangles (i.e. the image of 
${\cal S}(R_n^n)$ under the wiring of Fig.~\ref{wire}).
The set of monomials 
$(A_{i_1})^{j_1}(A_{i_2})^{j_2}\cdots(A_{i_p})^{j_p}$
where $i_k,\, j_k\in\NS$ and $\sum_{k=1}^pi_kj_k=n$
forms a linear basis for $\curl{C}^{(n)}$. (Note that
the number of such monomials is equal to the number of partitions of $n$.)
Thus, the algebra $\curlC^+$ is graded by weighted degree,
$\curlC^+=\bigoplus_{n=0}^\infty \curlC^{(n)}\,.$
\label{medc}

\subsection{Idempotents.}
\label{sec4}

In \cite{mine,idemp}, we 
adapt the construction of Gyoja \cite{gyoja}
(similar to that of Young symmetrisers) to
construct an idempotent element of 
$\curl{S}(R_n^n)$, for each Young diagram with $n$ cells.
Our building blocks are Jones idempotents \cite{hecke}
corresponding to single row and column Young diagrams. 
Equivalent idempotent elements appear in
the work of Yokota \cite{yoko}.  
For completeness, we give the bones of the construction.
For a full description, see \cite{mine,idemp}. 
Define $a_n$ and $b_n$ to be
\[
a_n=\sum_{\pi\in S_n} (-a)^{-l(\pi)}\omega_\pi\;,
\quad b_n=\sum_{\pi\in S_n} (-b)^{-l(\pi)}\omega_\pi\,,
\]
where $a=-xs^{-1}$ and $b=xs$ 
are the roots of the quadratic relation in $H_n$
and $l(\pi)$ is the length of the permutation $\pi$ (which
is equal to the writhe of the p.p.b. $\omega_\pi$).
In diagrams, we will represent the element
$a_l$ by a white box labelled $l$ and $b_k$ by a shaded box labelled $k$. 
A tensor product sign will denote the juxtaposition of two oriented tangles.
The elements $a_l$ and $b_k$ have the following properties.

\subsubsection{Lemma.{\rm\protect\cite{nato}}}
\label{linhom}

Let $\phi_a$ and $\phi_b$ be linear homomorphisms from $\curl{S}(R^n_n)$
to the ring of scalars $\Lambda$ defined by  
$\phi_a(\sigma_i)=a$ and 
$\phi_b(\sigma_i)=b$. Then for all $T\in \curl{S}(R^n_n)$,
\[
a_nT=Ta_n=\phi_b(T)a_n\,, \quad  b_nT=Tb_n=\phi_a(T)b_n\,.
\]
\qed
\subsubsection{Lemma.{\rm\protect\cite{mine}}}
\label{bill}

In $H_n$, we can decompose $a_l$ into a linear combination
of terms involving $a_{l-1}$:
\begin{eqnarray*}
a_l&=& a_{l-1}\otimes a_{1}\,+\,
                \sum_{i=0}^{l-2} (x^{-1}s)^{i+1}\ \
                a_{l-1}\sigma_{l-1}\sigma_{l-2}\cdots\sigma_{l-i-1}\\
        &=&\raisebox{-1cm}{\epsfysize.7in\epsffile{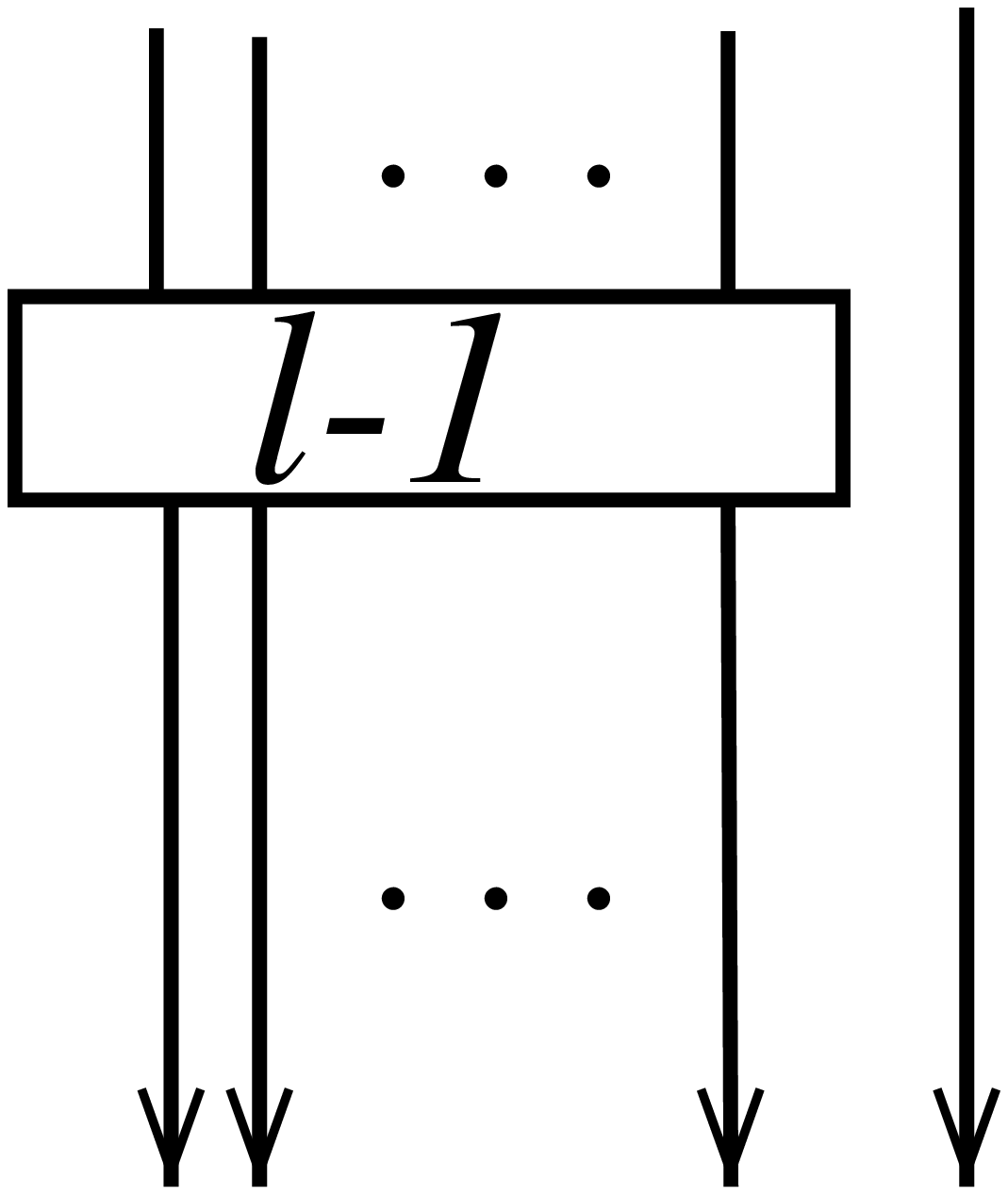}}\,+\,
                \sum_{i=0}^{l-2} (x^{-1}s)^{i+1}\ \
                \raisebox{-1cm}{\epsfysize.8in\epsffile{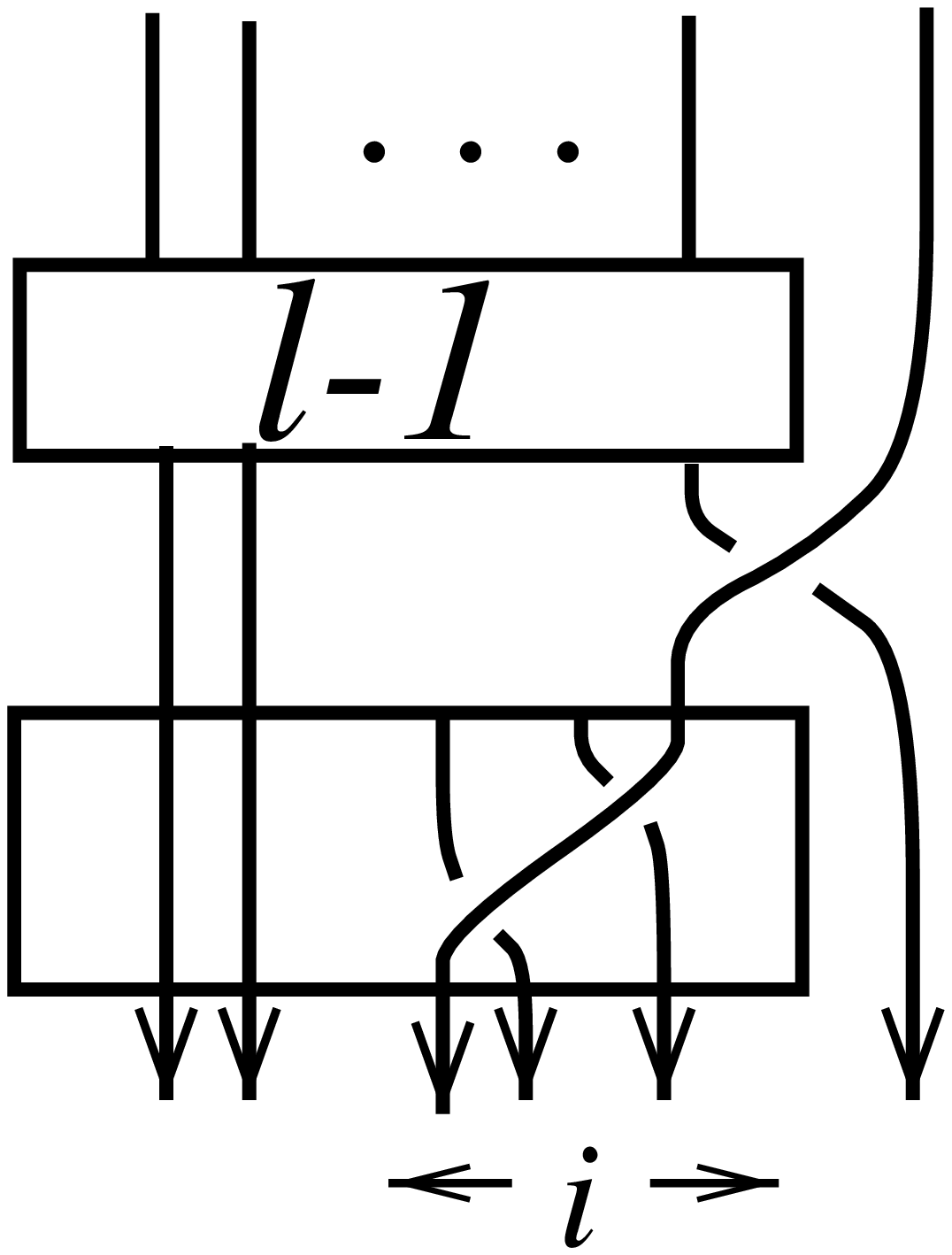}}\,.\\
\end{eqnarray*}
Similarly,
\begin{eqnarray*}
b_k&=& b_{k-1}\otimes b_{1}\,+\,
        \sum_{i=0}^{k-2}(-x^{-1}s^{-1})^{i+1}\ \
                b_{k-1}\sigma_{k-1}\sigma_{k-2}\cdots\sigma_{k-i-1}\\
        &=&\raisebox{-1cm}{\epsfysize.7in\epsffile{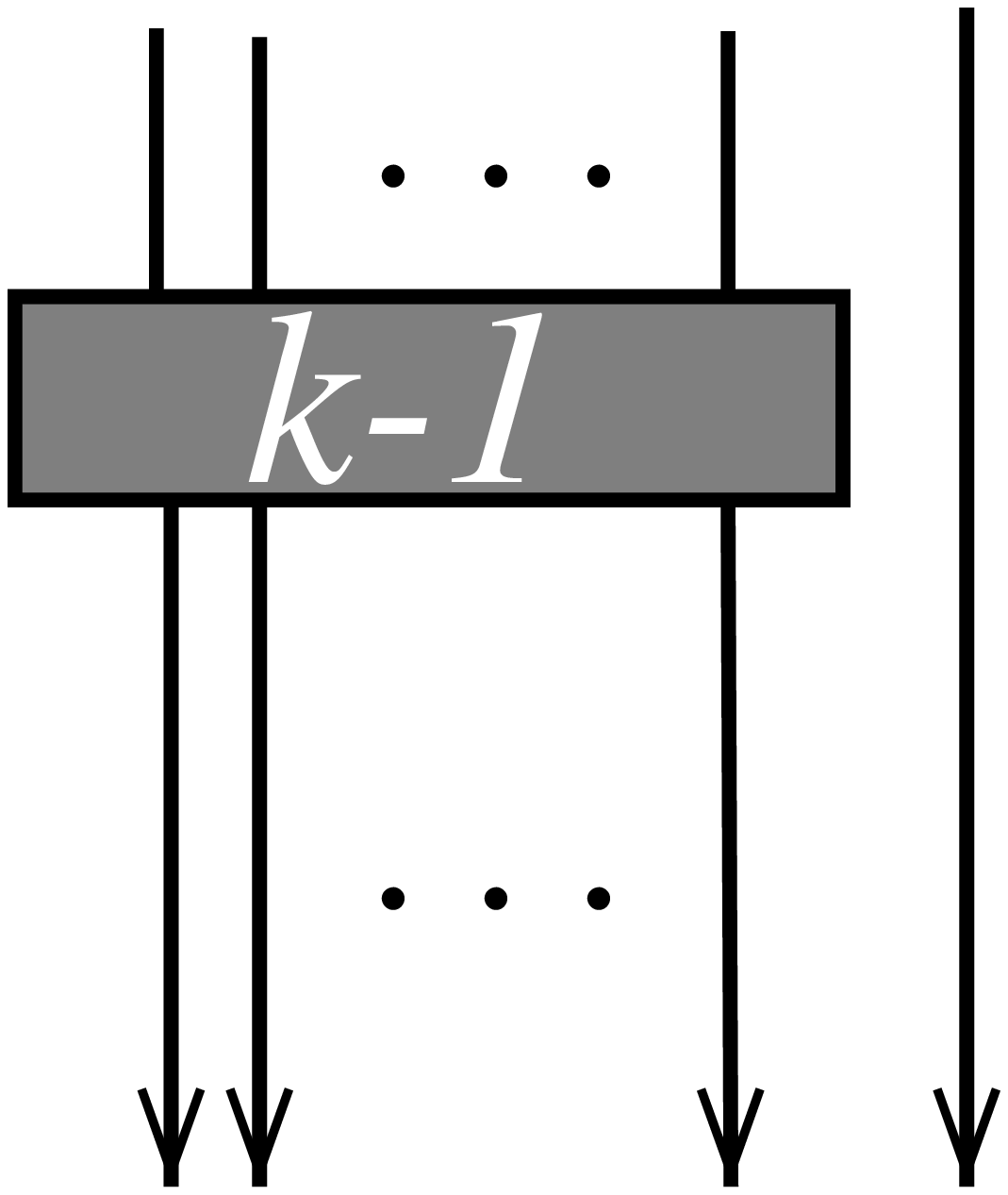}}\,+\,
        \sum_{i=0}^{k-2}(-x^{-1}s^{-1})^{i+1}\ \
                \raisebox{-1cm}{\epsfysize.8in\epsffile{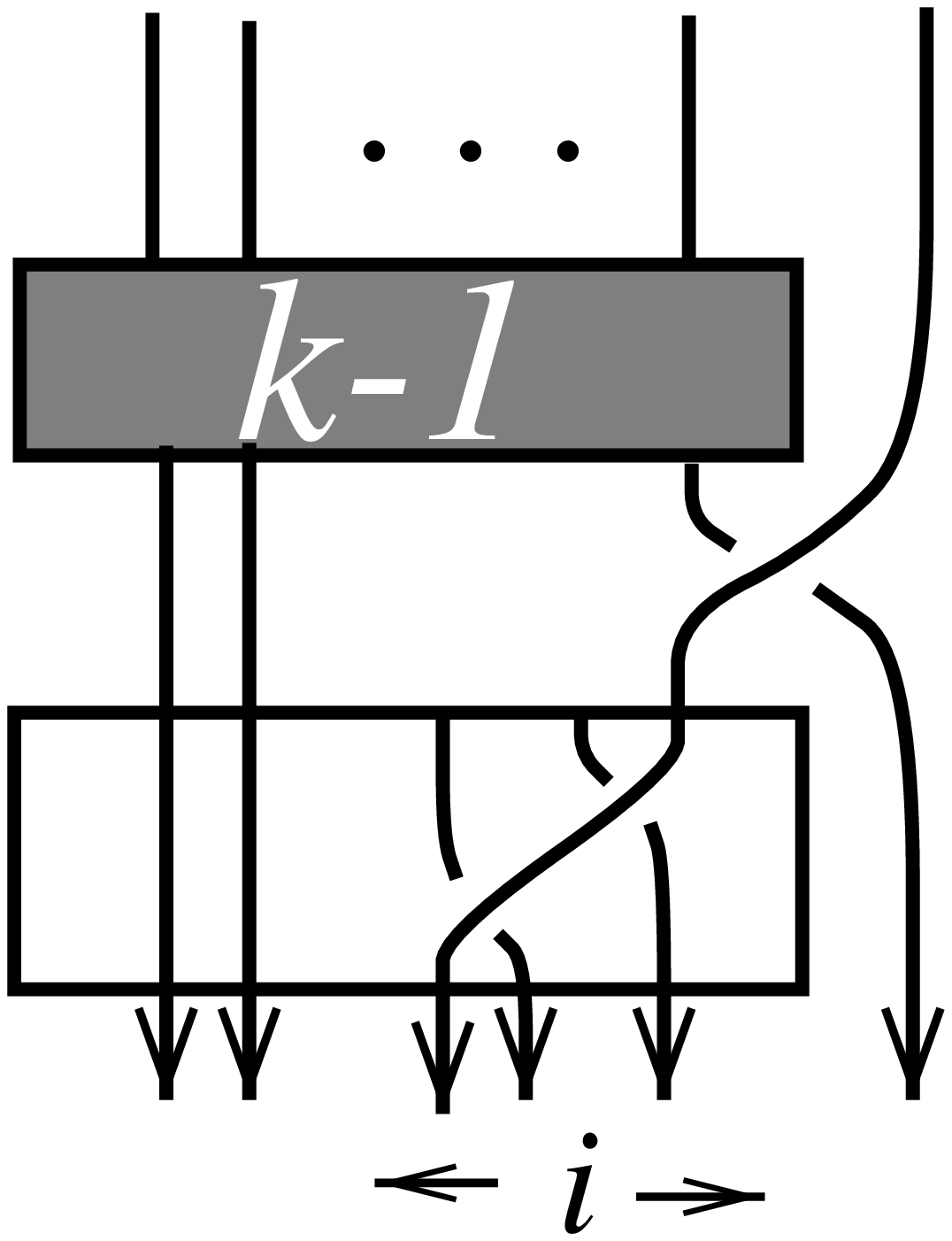}}\,.\\
\end{eqnarray*}
Equally, we can decompose $a_l$ and $b_k$ from above (turn the 
pictures upside down). 
\qed

\begin{figure}[ht]
\[
\epsfxsize1in\epsffile{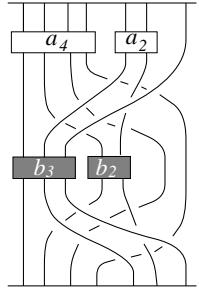}
\]
\caption{The quasi-idempotent $e_\nu$.}
\label{hats}
\end{figure}

We now define a quasi-idempotent element $e_\lambda$ for each $\lambda\in Y$. 
Assign a braid string to each cell, ordered by the tableau $T(\lambda)$.
On the strings in the $i$th row of $\lambda$, we place $a_{\lambda_i}$.
We will denote this linear combination of braids by $E_\lambda(a)$.
Similarly define $E_\lambda(b)$, replacing $a_{\lambda_i}$ by
$b_{\lambda_i}$.
For example, if $\nu=(4,2,1)$ then $E_\nu(a)=a_4\otimes a_2\otimes a_1$.
Define 
$e_\lambda=
E_\lambda(a)\omega_{\pi_\lambda}E_{\lambda^\vee}(b)\omega^{-1}_{\pi_\lambda}$,
where $\pi_\lambda$ is the permutation defined in Sect.~\ref{sec2}.
The element $e_\nu$ is shown in Fig.~\ref{hats}.

\subsubsection{Theorem.{\rm\protect\cite{mine,qdim}}}
\label{main}

The $e_\lambda$, $\vert\lambda\vert=n$,
are quasi-idempotent, mutually orthogonal elements of 
$\curl{S}(R^n_n)$ i.e. $e_\lambda^2=\alpha_\lambda\,e_\lambda$ for
some scalar $\alpha_\lambda$ and for $\lambda\neq\mu$, $e_\lambda e_\mu=0$.
We will denote the closures of the genuine idempotents, 
$1/\alpha_\lambda$ $e_\lambda$, 
by $Q_\lambda\in {\cal C}^+$.
The $Q_\lambda$ form a free $\Lambda$-basis for $\curl{C}^+$.
The value for $\alpha_\lambda$ is given by the formula
\[
\alpha_\lambda\ =\ \prod_{\mbox{\scriptsize{cells}}}
			s^{\mbox{\scriptsize{content}}}[\mbox{hook length}]
	\ =\ \prod_{(i,j)\in\lambda}s^{j-i}[\lambda_i+\lambda^\vee_j-i-j+1]
\]
Denote the $U_q(sl(N))$ invariant of a knot $K$ coloured by
the representation $V$ by $J(K;V)$. 
Let $C$ be a framed knot and let $S$ be the satellite $C*Q_\lambda$.
Set $\curl{X}_N$ to be the evaluation of $\cal X$ at
$x=s^{-1/N}$ and $v=s^{-N}$, where $s^2=q$.  Then,
\[
J(C;V_\lambda)=\curl{X}_N(S)\,.
\]
\qed

A subtle point (discussed in detail in \cite{knot96})
is that we require $1/\alpha_\lambda\in\Lambda$ and so
we must extend the scalar ring so that the quantum integers can be inverted.  
From now on, we shall assume that this is the case. 
As they will be used in the proceeding sections, we give the 
values of the scalar $\alpha$ for the hook shaped diagrams explicitly.
For these Young diagrams, it is quite simple to calculate $\alpha_{\mu_{k,l}}$
directly  using Lemma~\ref{bill}.
\[
\alpha_{\mu_{k,l}}
	=s^{(l(l-1)-k(k-1))/2}\,[\,k+l-1\,][\,k-1\,]![\,l-1\,]!\,.
\]

\subsection{An isomorphism from $\curlR_\infty$ to $\curl{C}^+$.}
\label{sec5}

We next describe an algebra isomorphism from
${\cal R}_\infty$ to $\curl{C}^+$.
We then demonstrate a generating set for
$\curl{C}^+$ by considering the image, under this isomorphism, of the 
Adams operators of $c_1$, which generate $\curl{R}_\infty$ as an algebra.

To improve notation marginally, we shall adopt $e_{k,l}$ to
denote the quasi-idempotent associated to the Young
diagram  $\mu_{k,l}$.  Similarly
we will write $Q_{k,l}$ and $\alpha_{k,l}$
rather than $Q_{\mu_{k,l}}$ and $\alpha_{\mu_{k,l}}$. 

\subsubsection{Lemma.}
\label{coeffak}

The element $Q_{k,1}$ is an element of ${\cal C}^{(k)}$ and can,
therefore, be expressed as a weighted homogeneous polynomial of degree $k$
in the $A_m$ with $k\geq m>0$. 
The coefficient of  $A_k$ in $Q_{k,1}$ is $(-x^{-1})^{k-1}/[k]!$.
\begin{sketchproof}

Since $Q_{k,1}$ is the closure of a linear combination of diagrams in
${\cal S}(R_k^k)$, $Q_{k,1}$ is an element of ${\cal C}^{(k)}$, and thus
weighted homogeneous of degree $k$, in terms of the $A_m$.

It can be shown by induction that $(-x^{-1})^{k-1}s^{-k(k-1)/2}[k-1]!$
is the coefficient of $A_{k}$ in $b_{k}$.
For the induction step, apply Lemma~\ref{bill} to $b_k$ to
obtain the term $b_{k-1}\otimes 1$ and a sum of $k$-string diagrams 
involving $b_{k-1}$.  The term $b_{k-1}\otimes 1$ has $A_1$ as a factor
thus, the monomial $A_k$ will not appear in this expression. 
Repeated application of Lemma~\ref{bill} to each of the remaining terms, 
along with some book keeping finishes the step.  Dividing by 
$\alpha_{k,1}$ we obtain  the coefficient of $A_k$ in $Q_{k,1}$ 
as $(-x^{-1})^{k-1}/[k]$.
\end{sketchproof}

\subsubsection{Proposition.}
\label{holger}

Define the algebra homomorphism 
$\theta:\curlR_\infty\rightarrow\curl{C}^+$
by $\theta(c_k)=Q_{k,1}\quad\forall k\in \NS$.
The homomorphism $\theta$ is an algebra isomorphism.
\begin{proof}
The ring $\curlR_\infty$ is defined to be the free
algebra generated by $c_k$, $k\in\NS$.
Hence, the homomorphism $\theta$ is an isomorphism if
we can show that the elements $Q_{k,1}$ freely generate
$\curl{C}^+$ as an algebra.

Recall that
$\{A_m : m\in\NS\}$, freely generate $\curl{C}^+$ as an algebra.
By induction on $n$, we will show that $A_n$
can be expressed in terms of the $Q_{k,1}$, for $0< k\leq n$.
and thus, $\{Q_{k,1}\ :\ k\in \NS\}$ generates $\curl{C}^+$ as an algebra.

For $n=1$, we have that $Q_{1,1}=A_1$, and we are done.

For $n>1$, assume that we have an expression for $A_m$
in terms of the $Q_{k,1}$ for all $m<n$.  
It follows from the definition of $\curl{C}^{(n)}$ that $Q_{n,1}\in\curl{C}^{(n)}$.
Therefore, there is an expression for $Q_{n,1}$ in 
terms of monomials of weighted degree $n$ in the $A_m$.
The monomials of weighted degree $n$ can be indexed by the partitions of
$n$, hence, by Lemma~\ref{coeffak}
\[
Q_{n,1}={(-x^{-1})^{n-1}\over[n]}\ A_n
	+\sum_{\vert\lambda\vert=n} 
		\beta_\lambda A_{\lambda_1}A_{\lambda_2}
						\cdots A_{\lambda_k}
\]
for some scalars $\beta_{\lambda}\in\Lambda$.
Thus, by the induction hypothesis, since the coefficient of $A_n$ is 
invertible in $\Lambda$, we can write $A_n$ in terms of the 
weighted degree $n$ monomials in the $Q_{k,1}$, $k\leq n$. Thus, these
monomials generate $\curl{C}^{(n)}$ as a $\Lambda$-module.  
This is equivalent to saying that the $Q_{k,1}$ generate $\curl{C}^+$ 
as an algebra.

It remains to show that $\curl{C}^+$ is generated freely by $Q_{k,1}$, 
$k\in\NS$.  We have just shown that the set of monomials
$Q_{i_1,1}^{j_1}Q_{i_2,1}^{j_2}\cdots Q_{i_m,1}^{j_m}$
for which $\sum_{p=1}^m i_pj_p=n$, generate $\curl{C}^{(n)}$
as a $\Lambda$-module.  By definition, $\curl{C}^{(n)}$
is freely generated by the set of weighted degree $n$ monomials in 
the $A_m$, which is a set with the same cardinality.
Therefore, since $\Lambda$ 
is commutative, the monomials in the $Q_{k,1}$ of weighted degree $n$
form a free $\Lambda$-basis for $\curl{C}^{(n)}$. 
It follows that $\{Q_{k,1}\,:\,k\in \NS\}$ generate the algebra 
$\curl{C}^+$ freely.
\end{proof}
\subsubsection{Lemma.}
\label{plum}

Let $\widehat{e}_{k,l}$ denote the closure of $e_{k,l}$ in $\curl{C}^+$.
The following relation holds :
\[
s^l [\,l\,] \widehat{e}_{k+1,l} + s^{-k} [\,k\,] \widehat{e}_{k,l+1}
	=s^{l-k}[\,l+k\,]\,\widehat{e}_{1,l}\,\widehat{e}_{k,1}\,.
\]
\begin{proof}
For simplicity, we draw the
oriented tangles which are wired into $\curl{C}^+$, as shown in Fig~\ref{wire}.
However, since we are working in $\curl{C}^+$, 
we can slide a piece of diagram off the top of a diagram and reintroduce it at the
bottom without loss.  All strings are assumed to be
oriented from top to bottom, unless indicated otherwise.
Applying Lemma~\ref{bill} to $\widehat{e}_{k+1,l}$,
\begin{eqnarray*}
\widehat{e}_{k+1,l}&=&\quad\raisebox{-1.3cm}{\epsfysize1in\epsffile{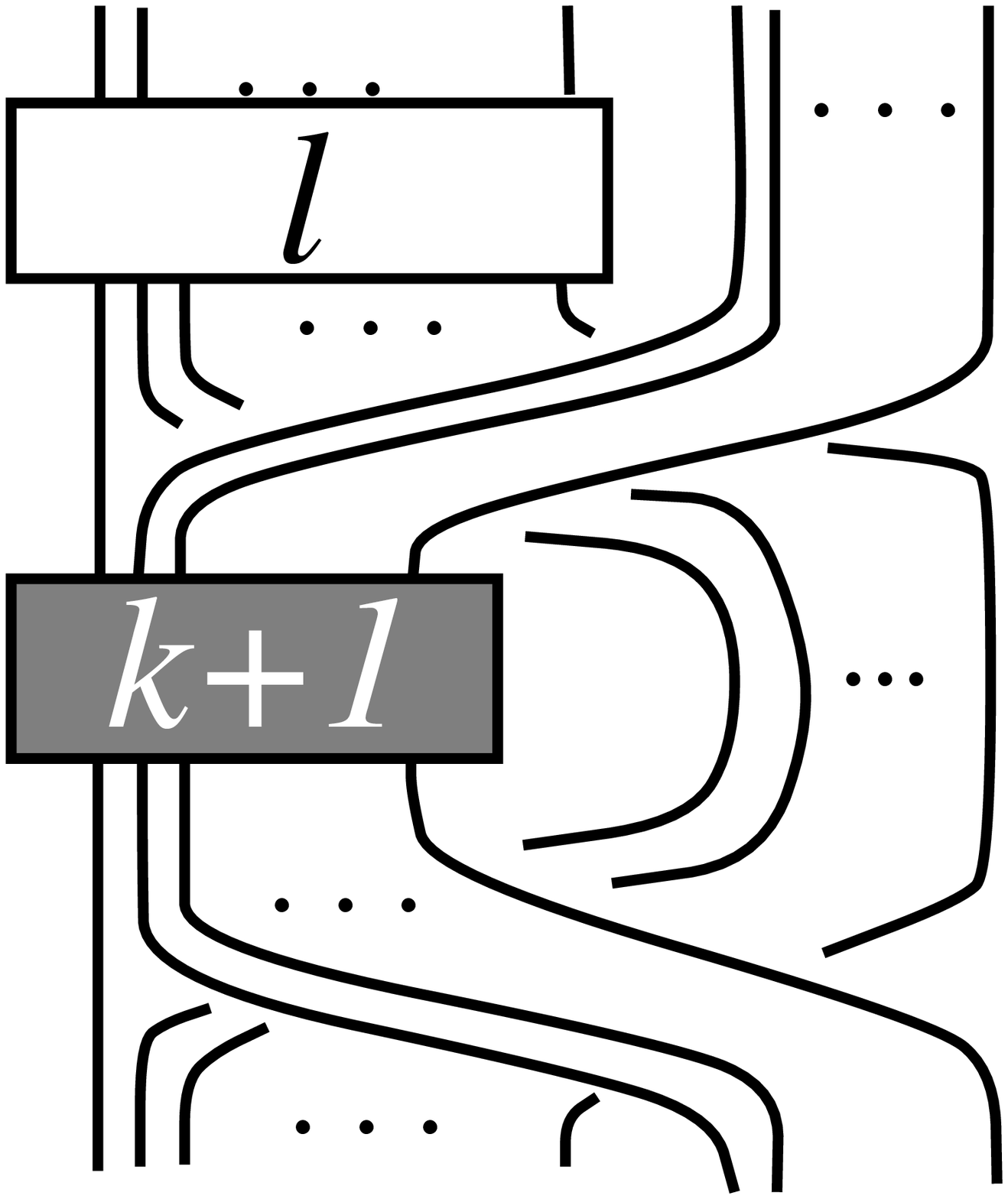}}
	 \quad=\quad\raisebox{-1.3cm}{\epsfysize1in\epsffile{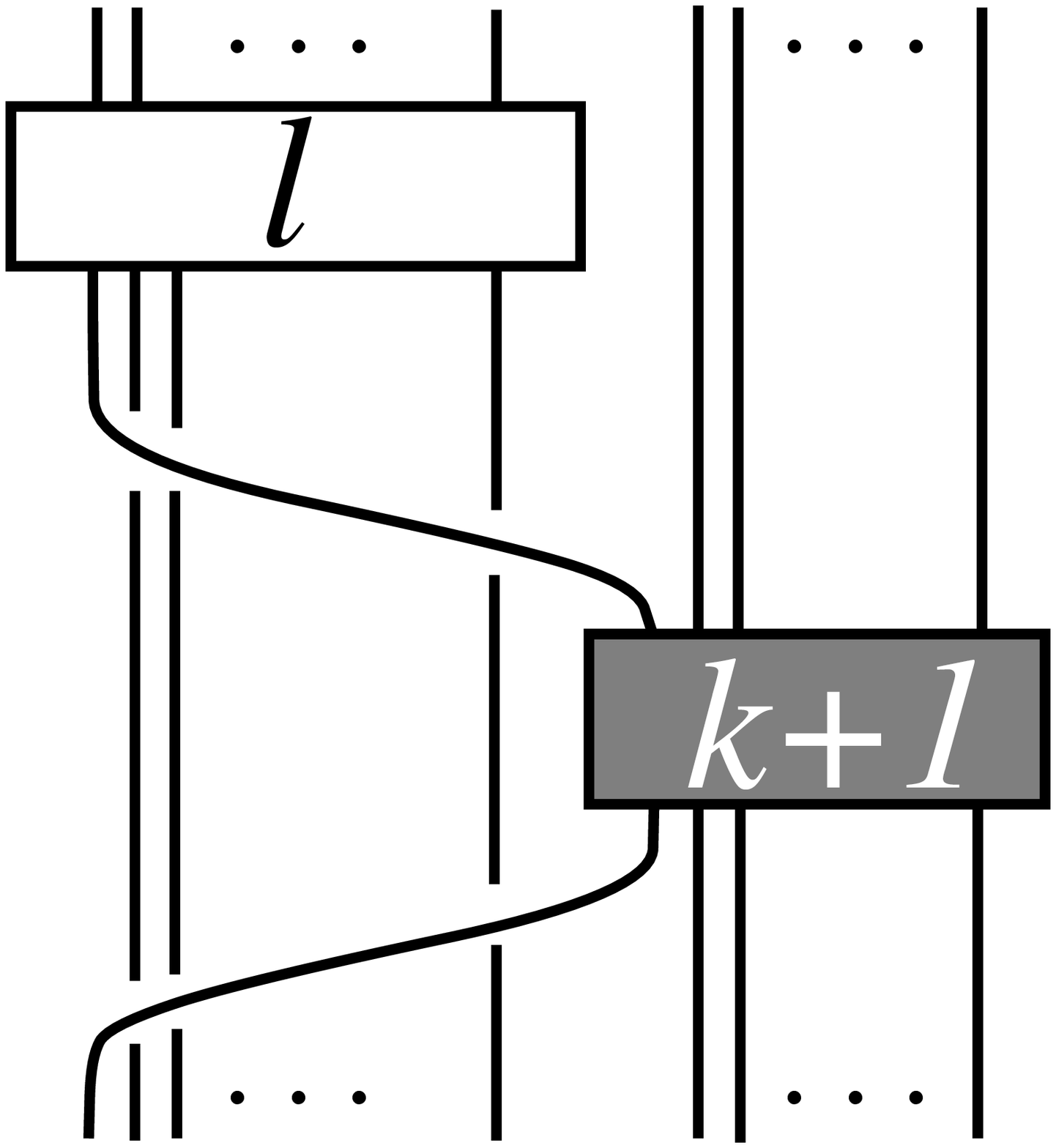}}\\
	 &=&\quad\raisebox{-1cm}{\epsfysize.9in\epsffile{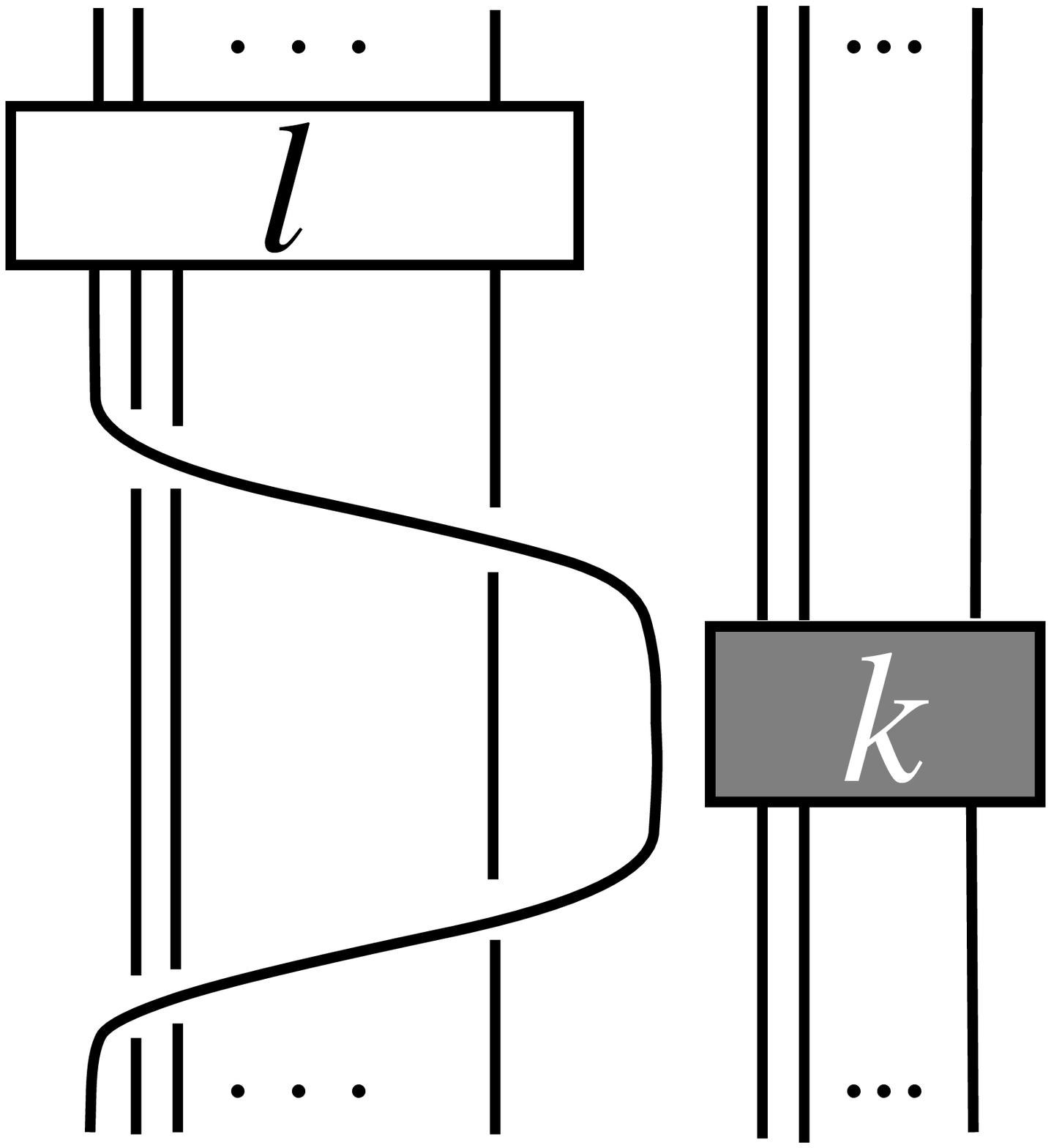}}
	 \quad+\quad\sum^{k-1}_{i=0}\left((-x^{-1}s^{-1})^{i+1}\ \raisebox{-1.3cm}
		{\epsfysize.9in\epsffile{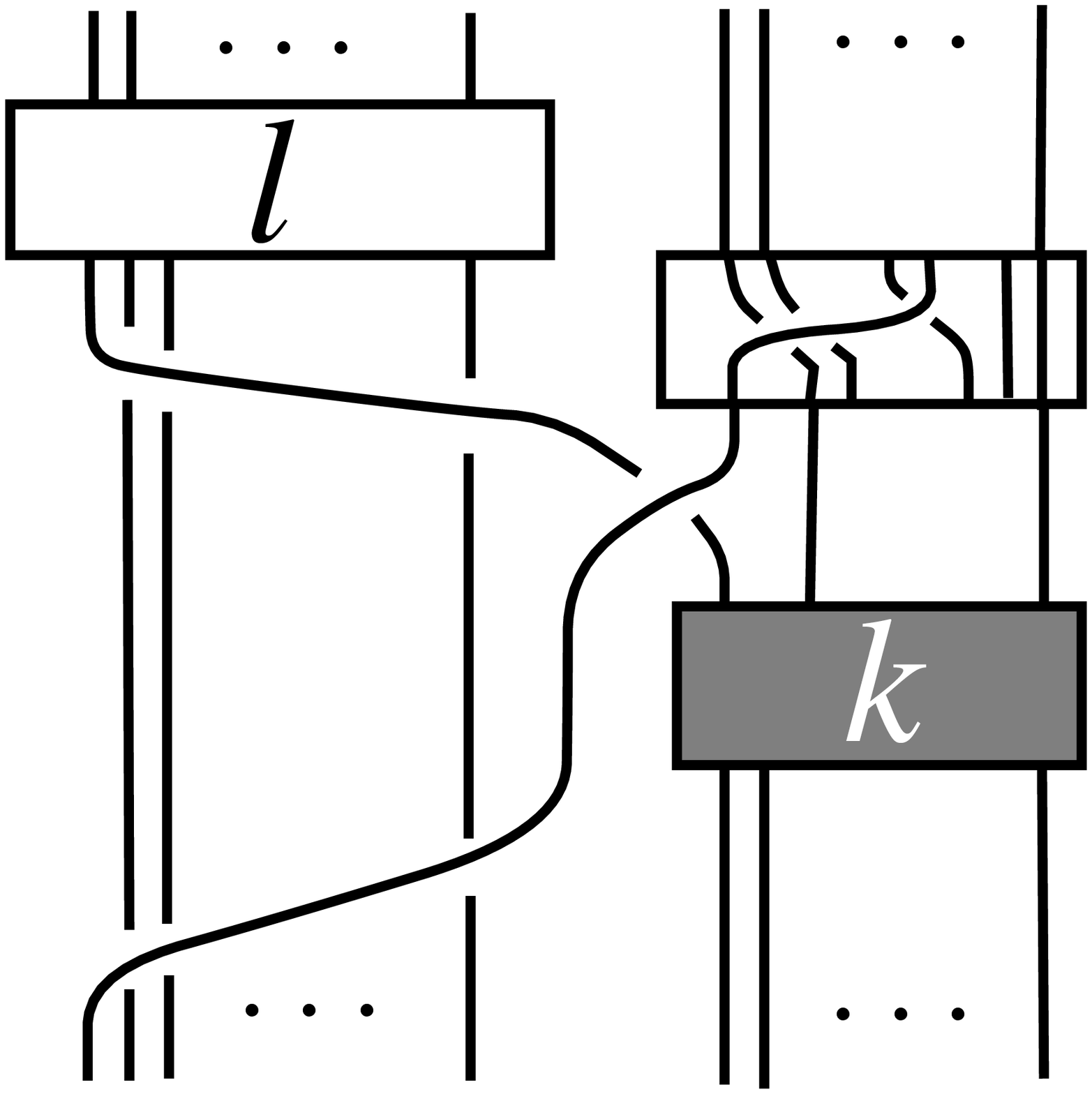}}\ \right)\,.\\
\end{eqnarray*}
Since we are working in the skein of the annulus, by Lemma~\ref{linhom},
\begin{eqnarray*}
\hat{e}_{k+1,l}&=&\raisebox{-1cm}{\epsfysize.9in\epsffile{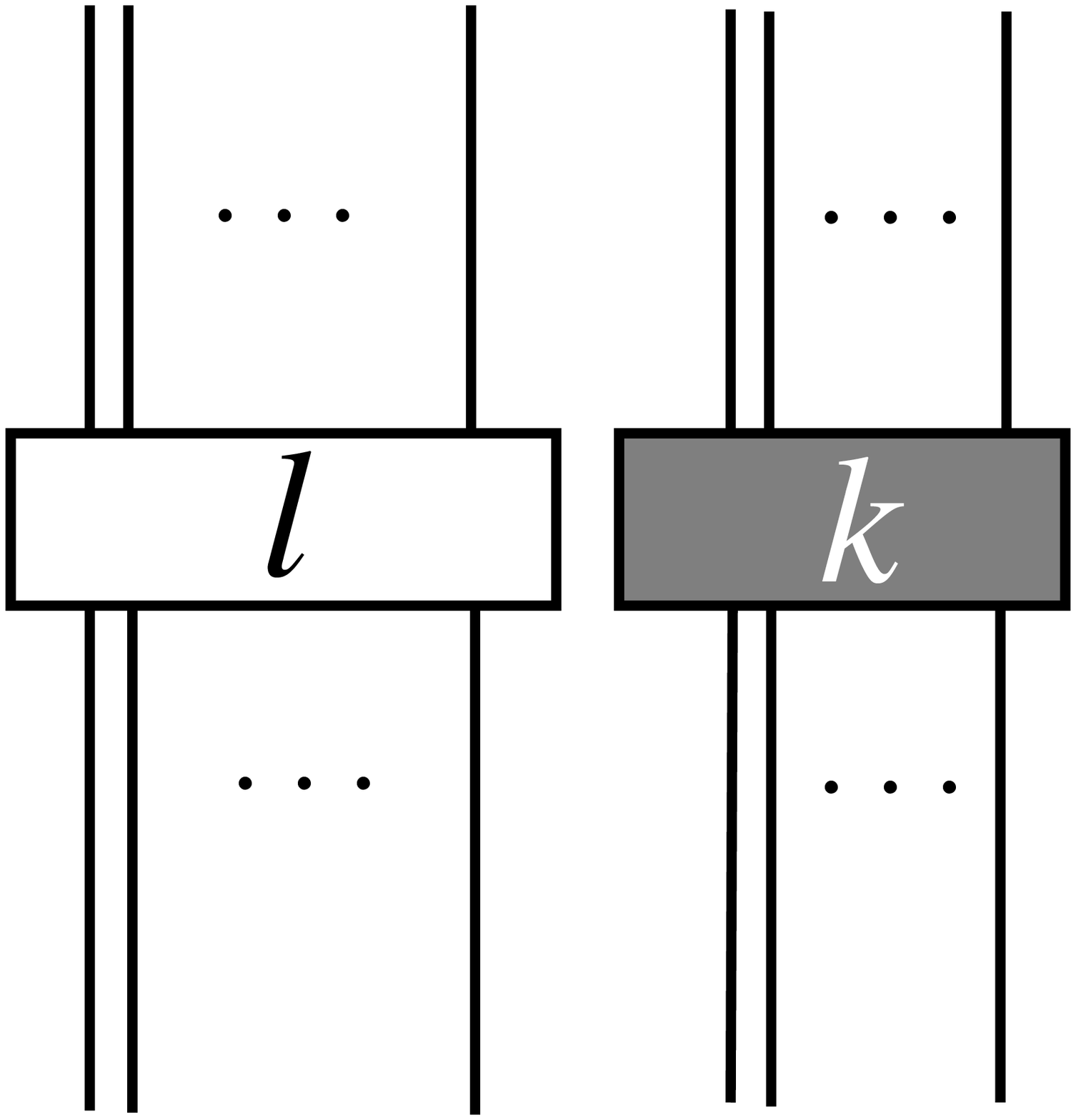}}\ 
   	     +\ \sum_{i=0}^{k-1}\left((-x^{-1}s^{-1})^{i+1}(-xs^{-1})^i\ 
		\raisebox{-1cm}{\epsfysize.9in\epsffile{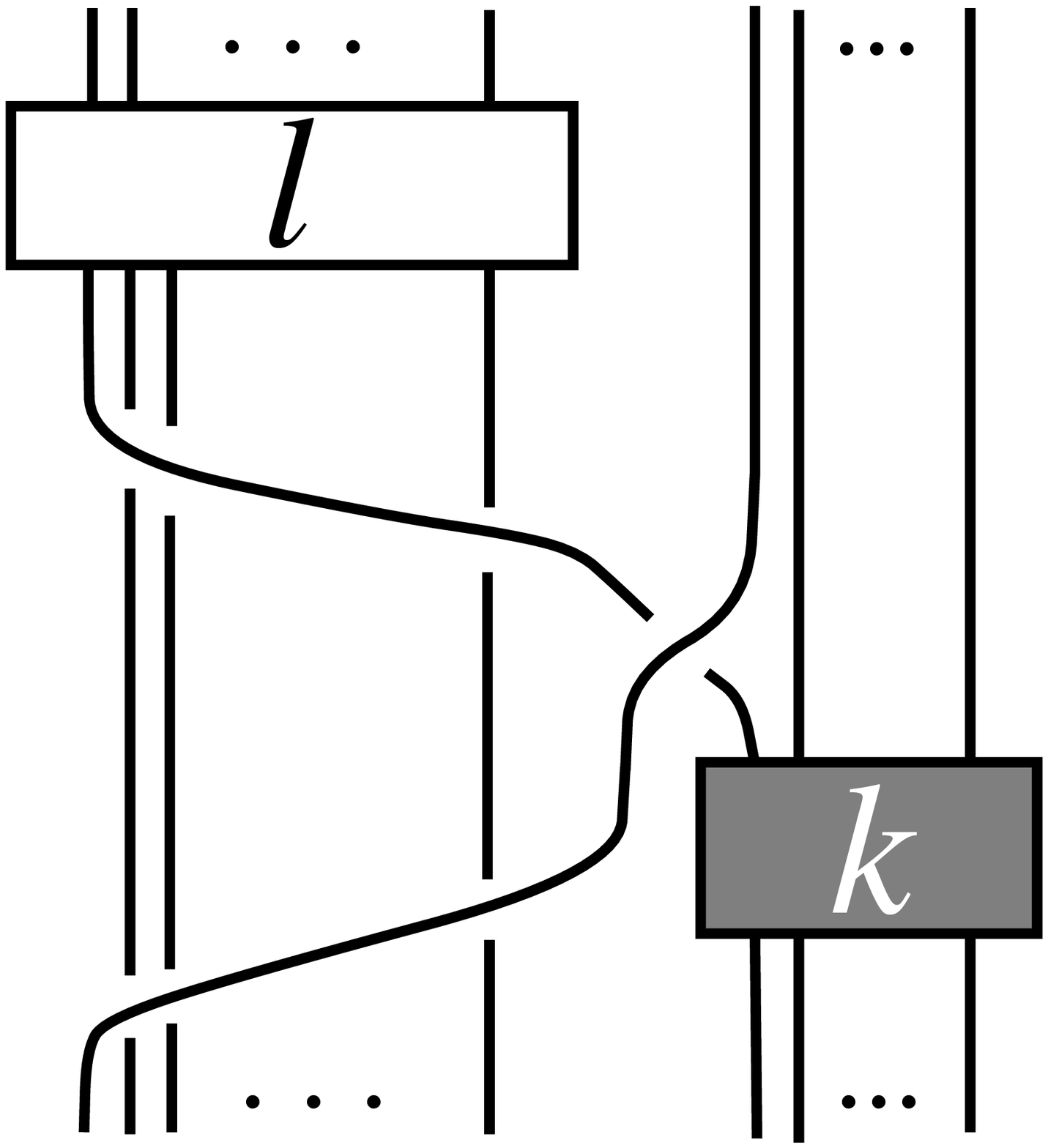}}\,
					\right)\\
	       &=&\raisebox{-1cm}{\epsfysize.9in\epsffile{emupard.ps}}\quad
		-\quad x^{-1}\left(\sum_{i=0}^{k-1}s^{-2i-1}\right)\quad
			\raisebox{-1cm}{\epsfysize.9in\epsffile{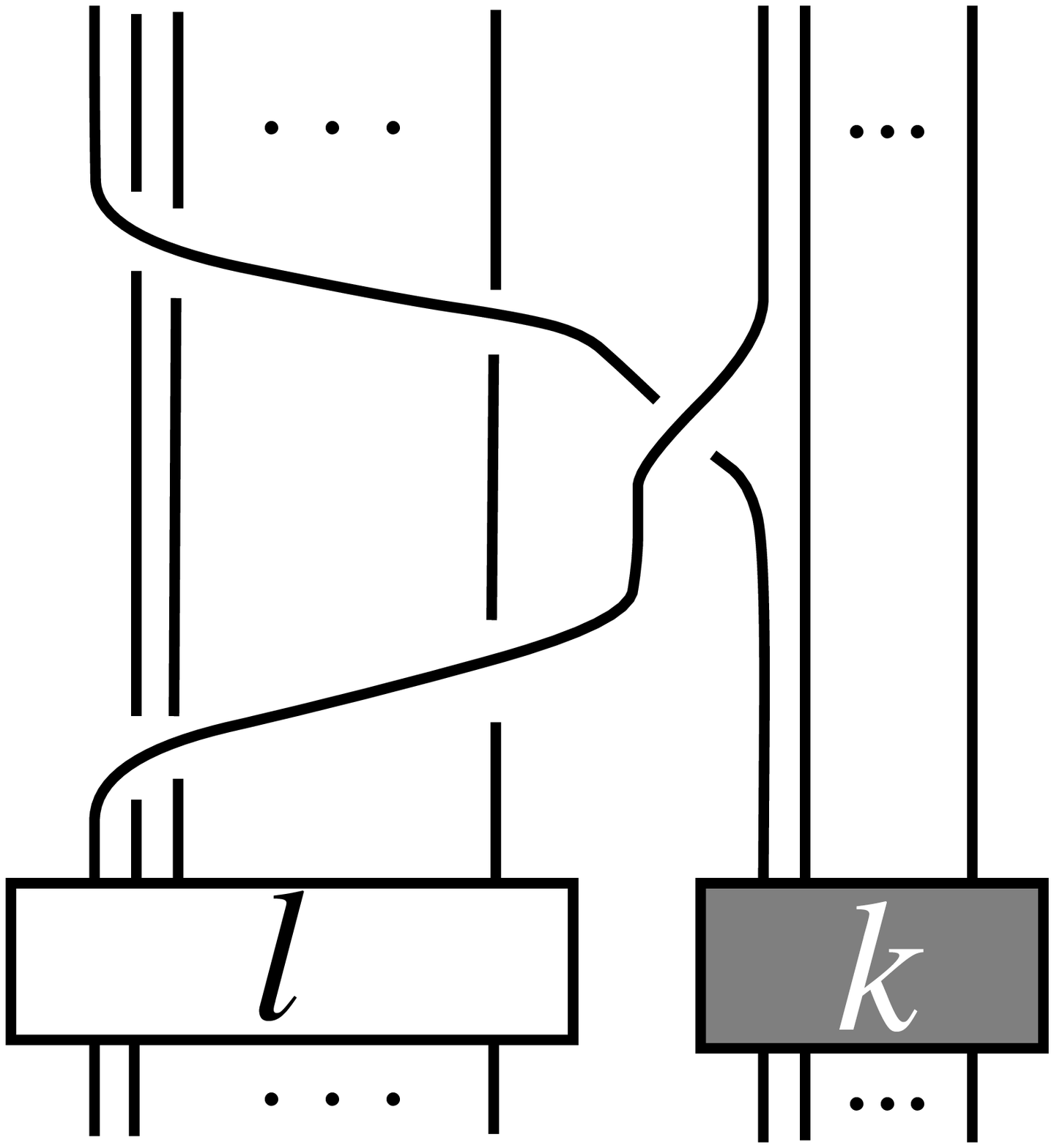}}\\
	       &=&\raisebox{-1cm}{\epsfysize.9in\epsffile{emupard.ps}}\quad
		-\quad x^{-1}s^{-k}[\,k\,]\quad
			\raisebox{-1cm}{\epsfysize.9in\epsffile{emubotd.ps}}\,.\\
\end{eqnarray*}
Using the same techniques on $\widehat{e}_{k,l+1}$ we obtain
\[
\widehat{e}_{k,l+1}\quad=\quad\raisebox{-1cm}{\epsfysize.9in\epsffile{emupard.ps}}\quad
			+\quad x^{-1}s^l[\,l\,] \quad
			\raisebox{-1cm}{\epsfysize.9in\epsffile{emubotd.ps}}\,.
\]
Hence, eliminating the term on the far right,
\[
s^l[\,l\,]\widehat{e}_{k+1,l}+s^{-k}[\,k\,]\,\widehat{e}_{k,l+1}
\ =\ (s^l[\,l\,]+s^{-k}[\,k\,])\,\widehat{e}_{1,l}\,\widehat{e}_{k,1}
\ =\ s^{l-k}\,[\,l+k\,]\,\widehat{e}_{1,l}\,\widehat{e}_{k,1}\,.
\]
\end{proof}
\subsubsection{Proposition.}
\label{bruce}

The idempotents $Q_{k,l}$ satisfy $Q_{k+1,l}+Q_{k,l+1}=Q_{k,1}Q_{1,l}$.

\begin{proof}
Substitute $\widehat{e}_{k,l}=\alpha_{k,l}\,Q_{k,l}$ into Lemma~\ref{plum}.
\end{proof}
\subsubsection{Corollary.}
\label{shep}

Let $Q_C(X)$ and $Q_D(X)$ respectively, denote the formal power series
\[
Q_C(X)=\sum_{k=0}^\infty (-1)^k Q_{k,1}X^k
\quad\mbox{ and }\quad
Q_D(X)=\sum_{l=0}^\infty Q_{1,l} X^l\,.
\]
(These formal power series can be thought of as
$\curl{C}^+$--versions
of the formal power series $C(X)$ and $D(X)$,
defined in Prop.~\ref{equal1}.)  Then $Q_C(X)$ is the inverse of
$Q_D(X)$,
\[
Q_C(X)Q_D(X)=1\,.
\]
\begin{proof}
Since Prop.~\ref{bruce} echoes Eq.~(\ref{ckdl}), 
the proof of Corollary \ref{shep} is analogous to that of Prop.~\ref{equal1}.
\end{proof}

We next demonstrate that $\theta(\mu_{k,l})=Q_{k,l}$.
(In fact it can be shown that for any Young diagram $\lambda$,
$\theta(\lambda)=Q_\lambda$, and that the $Q_\lambda$ form a linear
basis for $\curl{C}^+$, but we do not need this here. 
Details can be found in \cite{mine}.) 
First, we prove the following lemma.

\subsubsection{Lemma.}
\label{pixie}

The image of the Young diagram $(l)$ under $\theta$ is $Q_{1,l}$.
\begin{proof}
The proof goes by induction on the number of cells.
We know that $d_1=c_1$, therefore $\theta(d_1)=\theta(c_1)=Q_{1,1}$.

Assume that we have the result for $d_i$, for $i<m$.
By Prop.~\ref{equal1}, for $m\geq 1$,
$\sum_{k=0}^m (-1)^k c_kd_{m-k}=0$.
Since $\theta$ is an  algebra isomorphism, this implies that
\begin{equation}
\theta(d_m)=\sum_{k=1}^m (-1)^{k-1}\theta(c_k)\theta(d_{m-k})\,.
\label{thetadm}
\end{equation}
By Prop.~\ref{shep} and the induction hypothesis
\begin{eqnarray*}
0=\sum_{k=0}^m  (-1)^k Q_{k,1}Q_{1,m-k}
	&=&Q_{1,m}+ \sum_{k=1}^m (-1)^k \theta(c_k)\theta(d_{m-k})\\
	&=&Q_{1,m}-\theta(d_m)\ \ \mbox{ by Eq.~(\ref{thetadm}).}
\end{eqnarray*}
\end{proof}
\subsubsection{Proposition.}
\label{nows}

For $k,l\geq 1$, the image of $\mu_{k,l}$ under $\theta$, 
is $Q_{k,l}$.
\begin{proof}

By Lemma.~\ref{pixie}, $\theta(\mu_{1,l})=\theta(d_l)$ for $l\in\mathbb{N}$.
This provides the base for an induction on the number of cells
and the length of the first column.
Assuming the result for all hook shaped diagrams with fewer than $m$ cells 
or $m$ cells and at most $k$ cells in the first column, 
\begin{eqnarray*}
\theta\,(\mu_{k+1,m-k})&=&\theta(c_kd_{m-k})-\theta(\mu_{k,m-k+1})
				\quad\mbox{by Prop.~\ref{holger}}\\
	&=& Q_{k,1}\,Q_{1,m-k}-Q_{k,m-k+1}\quad\mbox{by the induction step}\\
	&=&Q_{k+1,m-k}\quad\mbox{by Prop~\ref{bruce}}\,,\\
\end{eqnarray*}
as required.
\end{proof}

\subsection{Adams operators.}
\label{sec6}

Recall that we can express the elements of $\curlR_N$ as
symmetric polynomials in $N$ variables.
It is well known that the power sums $\sum_i x_i^m$, $m\in\NS$,
generate the algebra of symmetric polynomials.
Thus, the images of these power sums under $\theta$, will give us an 
alternative generating set for ${\cal C}^+$.
We give an expression for the
$m$th power sum as a sum of $m$ closed braids in $\curl{C}^+$.
An advantage of these power sums over the $Q_{k,1}$ 
is that the number of braids in their expressions is linear 
rather than factorial in the number of strings and all the braids 
are reasonably simple.  
Recall the definitions of $C_N$ and $p_N$ from Sect.~\ref{sec2}.

\label{worm}
The \defn{Adams operators}, $\{\psi_m\}_{m\in \NS}$, 
are a family of  ${\cal R}_N$-endomorphisms, 
defined by their images on the $x_i$ (despite the fact that the $x_i$ are just
a formal device and are not elements of $\curl{R}_N$),
\[
\psi_m(x_i)=x_i^m.
\]
Hence $\psi_m(c_1)=\sum x_i^m$ is the $m$th Newton power sum.
As a polynomial in the $c_i$ it is independent of $N$, i.e. 
there is a polynomial $\beta_m({\bf c})\in \curlR_\infty$ 
with $p_N(\beta_m)=\psi_m(c_1)$ 
for all $N$.

\subsubsection{Proposition.}
\label{sum}

The following identities hold for $\psi_m(c_1)$,
\[
\psi_m(c_1)\,=\,p_N\left(\sum_{k=1}^m(-1)^{k-1}kc_kd_{m-k}\right)
	\,=\,p_N\left(\sum^m_{k=1}(-1)^{k-1}\mu_{k,m-k+1}\right)\,.
\]
\begin{proof}
The function $\ln(C(X))$ has formal power series expansion
\[
\ln(C(X))=\sum_{m=1}^\infty \beta_m({\bf c})X^m
\]
where $\beta_m({\bf c})$ denotes a polynomial in the $c_k$. 
Differentiating $\ln(C(X))$ with respect to $X$, we get 
\[
{C'(X)\over C(X)}= C'(X)D(X)=\sum_{m=1}^\infty m \beta_m({\bf c})X^{m-1}\,.
\]
By comparing coefficients it follows that 
\begin{eqnarray*}
m\beta_m({\bf c})&=&\sum_{k=1}^m (-1)^{k}kc_kd_{m-k}\\
   	 =&&\sum_{k=1}^{m-1}\left((-1)^{k}k(\mu_{k+1,m-k}+\mu_{k,m-k+1})\right) 
			\quad +(-1)^{m} m \mu_{m,1}\\
	=&&(-1)^{m} m \mu_{m,1}+\sum^{m-1}_{k=1}(-1)^{k}k\mu_{k,m-k+1} + 
		   \sum^m_{k=2}(-1)^{k-1}(k-1)\mu_{k,m-k+1}  \\
		 =&&-\mu_{1,m}+(-1)^{m}\mu_{m,1}+
			\sum_{k=2}^{m-1}(-1)^k\mu_{k,m-k+1} \\
		 =&&\sum^m_{k=1}(-1)^{k}\mu_{k,m-k+1} \quad.
\end{eqnarray*}
For each $N$ we have that
\[
\ln(C_N(X))   \ =\ \sum_{i=1}^N\ln(1-x_iX)
		\ =\ -\sum_{m= 1}^\infty\sum_{i=1}^N{x_i^m\over m} X^m
		\ =\ - \sum_{m=1}^{\infty}{\psi_m({c_1})\over m}X^m\,.
\]
Now $\ln(C_N(X))=p_N(\ln(C(X)))$, therefore, 
$\psi_m(c_1)=- p_N(m \beta_m({\bf c}))$, as required.
\end{proof}
Since these formulae are independent of $N$, 
for most purposes we can treat $\psi_m(c_1)$ as if 
it were an element of $\curlR_\infty$.

\subsubsection{Lemma.}
\label{rich}
Let $\Psi(X)=\sum_{m=1}^\infty\psi_m(c_1)X^{m-1}$.
Then 
\[
\Psi(X)=-C\/'(X)D(X)=D\/'(X)C(X)\,.
\]
\begin{proof}
From Prop.~\ref{sum} we know that 
\[
\Psi(X)=-\,{d\ \over dX}\,\ln(C(X))=-C\/'(X)D(X)\,.
\]
It remains to note that, by Prop.~\ref{equal1}
\[
D\/'(X)C(X)=\left({1\over C(X)}\right)^{'}C(X)
		={-C\/'(X)\over C(X)^2}C(X)=-C\/'(X)D(X)\,.
\]
\end{proof}

Strickland and independently Rosso and Jones \cite{strick,rojo} 
proved the following formula for the $U_q(sl(N))$-invariants of torus knots.
However, the independence of the result from the knot $K$ implies
that the relation also holds directly in ${\cal R}_\infty$ (and via $\theta$
in $\curl{C}^+$).

\subsubsection{Theorem. {\rm \protect\cite{strick,rojo}}}
\label{206}

Let $m$ and $p$ be coprime integers.  Let $K^{(m,p)}$ denote the $(m,p)$
cable of the knot $K$.  Then if 
$\psi_m(c_1)=\sum_{\tau\in Y}\beta_\tau \tau$ then
\[
J(K^{(m,p)};c_1)=\sum_\tau \beta_\tau f_\tau^{p/m}J(K;\tau)\,.
\]
where the $(m,p)$-cable of $K$ is the satellite of $K$ with pattern
the $p$th power of the tangle
\[
\begin{array}{c}
\epsfysize.5in\epsffile{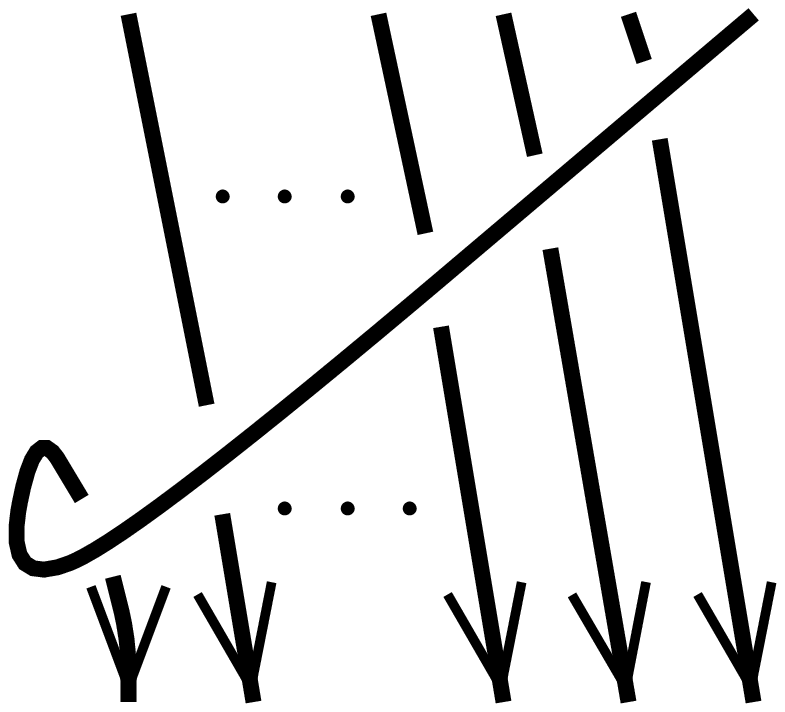}\\
\leftarrow\ m\ \rightarrow\\
\end{array}
\]
and $f_\tau$ is the framing factor associate to $\tau$ and $J(K;\tau)$
is the $U_q(sl(N))$ invariant of $K$ coloured by $\tau$.
\qed

By Prop.~\ref{sum} and Theorem~\ref{206}, setting $p=1$, we obtain
the following corollary.

\subsubsection{Corollary.}
\label{sparks}
\[
A_m =
(xv^{-1})^{-1}\sum_{k=1}^m (-1)^{k-1}f_{k,m-k+1}^{{1\over m}}
					\theta(\mu_{k,m-k+1})\,,
\]
where $f_{k,m-k+1}$ is the framing factor of the Young diagram $\mu_{k,m-k+1}$.
\qed

Originally calculated in \cite{morton}, the framing factors for $U_q(sl(N))$ 
calculated via skein theory \cite{idemp} are given by  
\[
f_{\lambda} = x^{\vert\lambda\vert^2}v^{-\vert\lambda\vert}s^{n_\lambda}\,,
\]
where $n_\lambda=\sum_i\lambda_i^2-\sum_j(\lambda_j^\vee)^2$.
For $\mu_{k,m-k+1}$, $n_{k,m-k+1}=m^2-2mk+m$.  Note 
that, since $\mu_{k,m-k+1}$ has $m$ cells, 
$f_{k,m-k+1}^{{1\over m}}$ contains only integer powers of 
$x$, $v$ and $s$.  We can obtain a similar expression for the 
closure of the $m$-string n.p.b., $ \overline{A}_{m}$,
by replacing $s$, $x$ and $v$ by $s^{-1}$, $x^{-1}$ and $v^{-1}$ respectively 
in the formula for $ A_m$. 
To see this substitute $p=-1$ in Theorem~\ref{206} and note that in 
$\curl{C^+}$ the closure of $\sigma_1\sigma_2\cdots\sigma_{m-1}$ is equivalent
to the closure of $\sigma_{m-1}\sigma_{m-2}\cdots\sigma_{1}$.

\subsubsection{Proposition.}    
\label{karen}

The following equalities hold:
\begin{eqnarray}
 A_m&=&x^{m-1}\sum^{m}_{k=1}(-1)^{k-1}s^{m-k}[k]\theta(c_kd_{m-k})\,,
\label{ameq1}
\\
 \overline{A}_{m}&=&x^{-(m-1)}\sum^{m-1}_{k=0}(-1)^{k}s^k[m-k]\theta(c_kd_{m-k})\,.
\label{ameq2}
\end{eqnarray}
\begin{proof}
By Cor.~\ref{sparks}
\begin{eqnarray*}
 A_m
        &=&x^{-1}v\sum_{k=1}^m(-1)^{k-1}x^mv^{-1}s^{m-2k+1}
                                                \theta(\mu_{k,m-k+1})\\
        &=&\sum^m_{k=1}(-1)^{k-1}x^{m-1}s^{m-2k+1}
                                                \theta(\mu_{k,m-k+1})\\
        &=&x^{m-1}\sum^m_{k=1}(-1)^{k-1}s^{m-2k+1}\sum_{i=k}^m(-1)^{i-k}
                                                \theta(c_id_{m-i})\\
        &=&x^{m-1}\sum^m_{i=1}\sum_{k=1}^i(-1)^{i-1}s^{m-2k+1}
                                                \theta(c_id_{m-i})\\
    &=&x^{m-1}\sum_{i=1}^m(-1)^{i-1}s^{m+1}\left(\sum_{k=1}^i s^{-2k}\right)
                                                \theta(c_id_{m-i})\\
    &=&x^{m-1}\sum^m_{i=1}(-1)^{i-1}s^{m+1}s^{-2}{(1-s^{-2i})\over (1-s^{-2})}
                                                \theta(c_id_{m-i})\\
        &=&x^{m-1}\sum^{m}_{i=1}(-1)^{i-1}s^{m-i}[i]\theta(c_id_{m-i})\,.\\
\end{eqnarray*}
This establishes Eq.~\ref{ameq1}. For Eq.~\ref{ameq2}, note that
$\,\sum^m_{i=0}(-1)^i\theta(c_i)\theta(d_{m-i})=0\,$ by Prop.~\ref{pixie}, 
therefore,
\begin{eqnarray*}
 A_m&=& A_m\, +\,\theta(x^{m-1}[m]\sum_{i=0}^m (-1)^ic_id_{m-i})\\
        &=&[m]x^{m-1}\theta(c_0d_{m})+x^{m-1}\sum_{i=1}^m(-1)^{i-1}
                                        (s^{m-i}[i]-[m])\theta(c_id_{m-i})\\
        &=&[m]x^{m-1}\theta(c_0d_m)+x^{m-1}\sum_{i=1}^m(-1)^{i}\left(
                     \sum_{k=1}^{m-i} s^{-m-1+2k}\right)\theta(c_id_{m-i})\,.\\
\end{eqnarray*}
Then
\begin{eqnarray*}
 \overline{A}_{m}&=&x^{-(m-1)}\sum^m_{i=0}(-1)^{i}\left(\sum^{m-i}_{k=1}
                                        s^{m+1-2k}\right)\theta(c_id_{m-i})\\
           &=&x^{-(m-1)}\sum^m_{i=0}(-1)^{i}s^i[m-i]\theta(c_id_{m-i})\\
           &=&x^{-(m-1)}\sum^{m-1}_{i=0}(-1)^{i}s^i[m-i]\theta(c_id_{m-i}).\\
\end{eqnarray*}
\end{proof}

We define $\Phi^+(X)$ and $\Phi^-(X)$ to be the formal power series
\begin{equation}
\label{lucy}
\Phi^+(X)=\sum_{m=1}^{\infty}  A_m X^{m-1}\quad\mbox{ and }\quad
\Phi^-(X)=\sum_{m=1}^{\infty} \overline{A}_{m} X^{m-1}\quad.
\end{equation}
We also define ``quantum derivatives'' of $C(X)$ and $D(X)$,
\[
{C'}_q (X)=\sum_{k=1}^{\infty}(-1)^{k}[k]c_k\,X^{k-1} \mbox{\quad and \quad}
                        {D'}_q(X)=\sum_{l=1}^{\infty}\,[l\/]d_l\,X^{l-1}\ .
\]
\subsubsection{Proposition.}
\label{mary}

The following three identities hold:
\begin{eqnarray*}
\Phi^+(X)&=&\,-\,\theta\left(\,{C'}_q(xX)D(xsX)\,\right)\,,\\
\Phi^-(X)&=&\theta\left(\,C(x^{-1}sX){D'}_q(x^{-1}X)\,\right)\,,\\
\Phi^-(X)&=&\,-\,\theta\left(\,{C'}_q (x^{-1}X)D(x^{-1}s^{-1}X)\,\right)\,.
\end{eqnarray*}
\begin{proof}
The first two identities follow automatically from
Prop.~\ref{karen}.
To see the third identity note that, by 
Prop.~\ref{equal1},
\[
\sum^\infty_{m=1}x^{-m+1}[m]
                        \left(\sum^m_{k=0} (-1)^k c_kd_{m-k}\right)
                                        X^{m-1}\quad=\quad0\,.
\]
We can, therefore, subtract any multiple of the image of this sum under
$\theta$ from $\Phi^-(X)$.  Hence
\begin{eqnarray*}
\Phi^-(X)&=&\Phi^-(X)-\sum_{m=1}^\infty x^{-m+1}\sum_{k=0}^m
                   (-1)^k[m]\theta(c_kd_{m-k})\,X^{m-1}\\
     &=&\sum_{m=1}^\infty x^{-m+1}\sum_{k=0}^m
                   (-1)^k(s^k[m-k]-[m])\theta(c_kd_{m-k})\,X^{m-1}\\
     &=&\sum_{m=1}^\infty x^{-m+1}\sum_{k=0}^m
              (-1)^k\left({s^m-s^{-m+2k}-s^m+s^{-m}\over s-s^{-1}}\right)
                                                \theta(c_kd_{m-k})\,X^{m-1}\\
     &=&\sum_{m=1}^\infty x^{-m+1}\sum_{k=0}^m
                   (-1)^ks^{-m+k}\left({-s^k+s^{-k}\over s-s^{-1}}\right)
                                                \theta(c_kd_{m-k})\,X^{m-1}\\
     &=&\sum_{m-1}^\infty x^{-m+1}\sum_{k=0}^m
                   (-1)^{k-1}s^{-m+k}[k]\theta(c_kd_{m-k})X^{m-1}\,.\\
\end{eqnarray*}
The third equation now follows by comparing the coefficient of $X^{m-1}$
above with that of $X^{m-1}$ in 
$\theta\left(\,C(x^{-1}sX){D'}_q(x^{-1}X)\,\right)$.
\end{proof}

Let $A_{i,j}$ be the closure of the braid below in $\curl{C}^+$,
with $i$ positive crossings and $j$ negative crossings
(so $A_m= A_{m-1,0}$ and $\overline{A}_{m}=A_{0,m-1}$).
\[
A_{i,j}\ =\ \raisebox{-8mm}{\epsfxsize.5in\epsffile{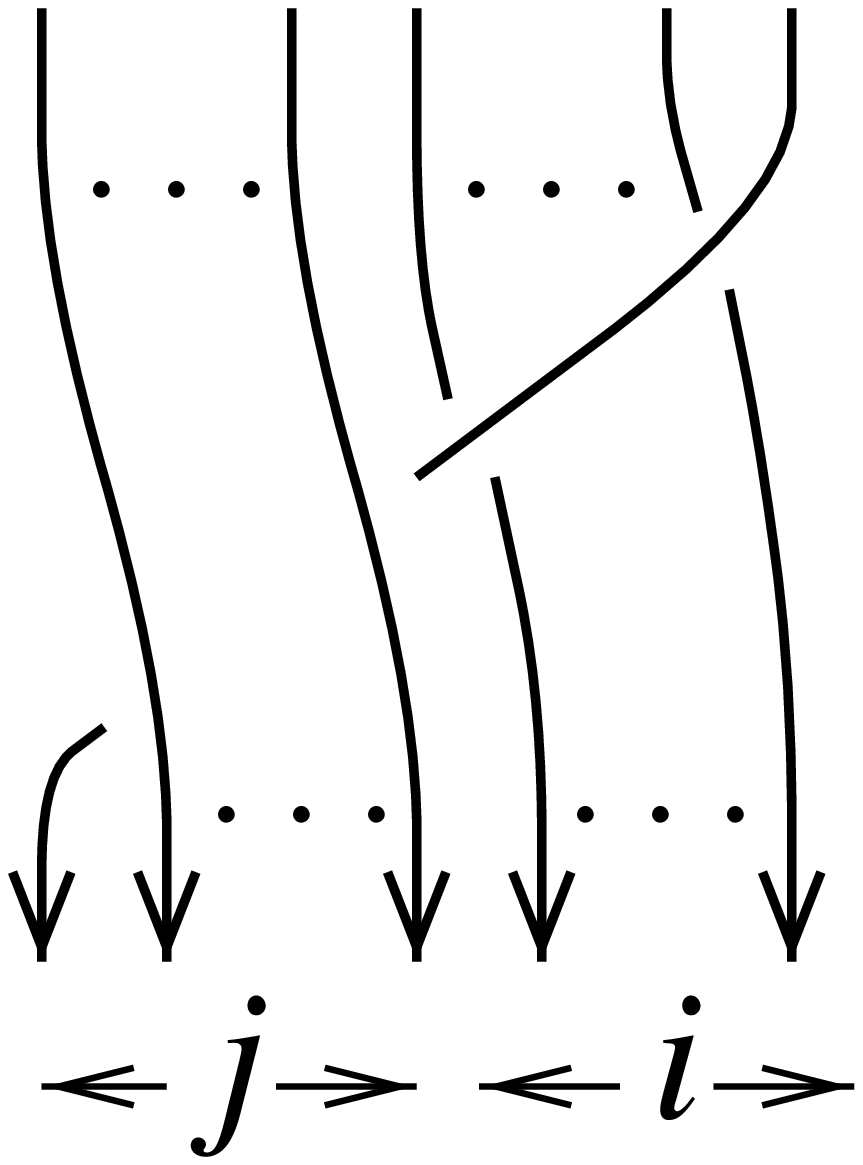}}
\]
\subsubsection{Lemma.}
\label{mark}
Let 
\begin{eqnarray*}
P_m&=&x^{m-1}\,\raisebox{-2mm}{\epsfxsize0.3in\epsffile{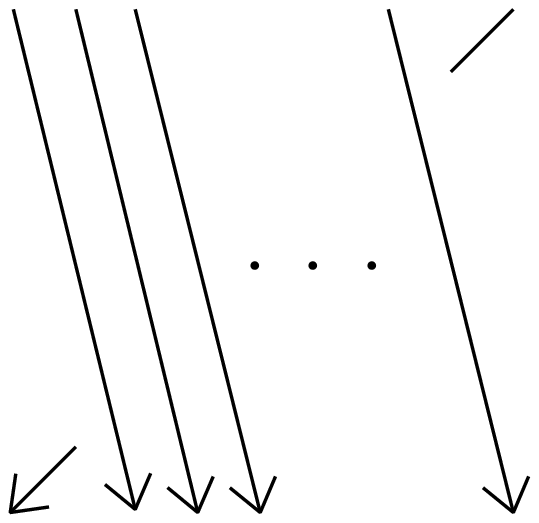}}
		   \,+\,x^{m-3}\,\raisebox{-2mm}
				{\epsfxsize0.3in\epsffile{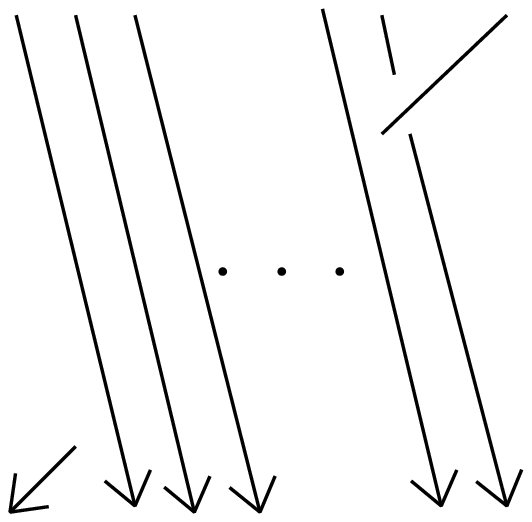}}
		   \,+\,x^{m-5}\,\raisebox{-2mm}
				{\epsfxsize0.3in\epsffile{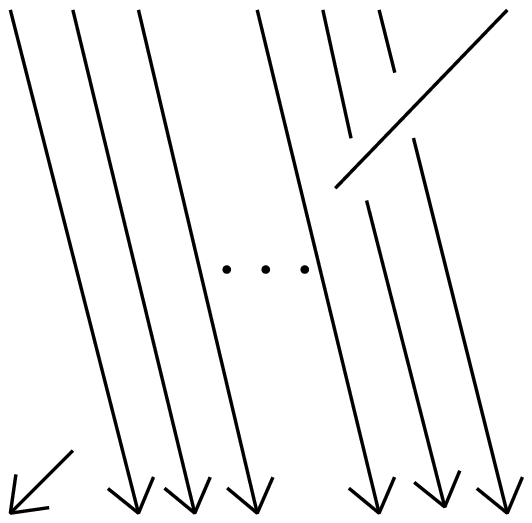}}
		   \,+\,\cdots\,
		     +\, x^{-m+3}\,\raisebox{-2mm}
				{\epsfxsize0.3in\epsffile{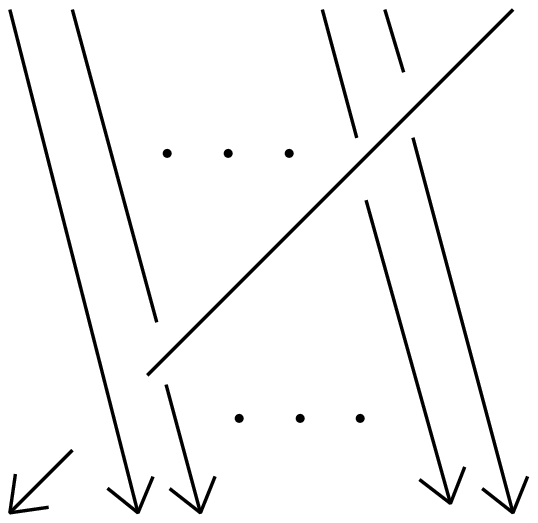}}
		   \,+\,x^{-m+1}\raisebox{-2mm}
				{\epsfxsize0.3in\epsffile{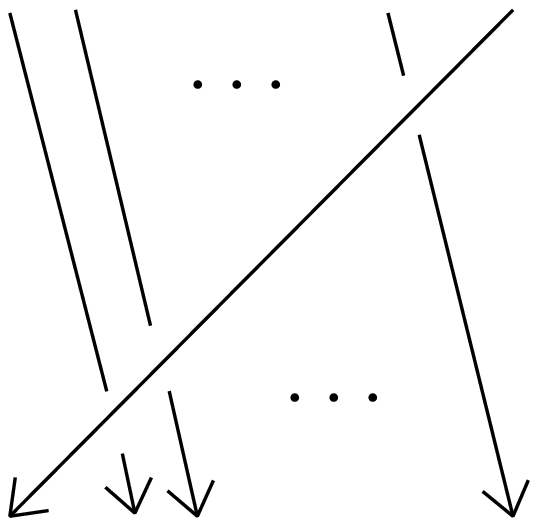}}\\
&=&\sum_{i=0}^{m-1}x^{m-1-2i} A_{i,m-1-i}\ .
\end{eqnarray*}
Then 
$P_m=\displaystyle\sum_{k=1}^{m-1}
		(s-s^{-1})x^{m-1-k}P_k A_{0,m-1-k}\ +mx^{m-1} A_{0,m-1}$.
\begin{proof}
The method here is to switch and smooth the positive crossings in 
the braids using the
skein relation. Numbering the crossings from bottom left to top right on 
the diagrams, we  start with the $(m-1)$st (positive) crossing.
\begin{eqnarray*}
P_m	&=& 	   x^{m-1}\raisebox{-2mm}
				{\epsfxsize0.3in\epsffile{basic1.ps}}
		   \,+\, x^{m-3}\raisebox{-2mm}
				{\epsfxsize0.3in\epsffile{basic2.ps}}
		   \,+\, x^{m-5}\raisebox{-2mm}
				{\epsfxsize0.3in\epsffile{basic3.ps}}
		   \,+\,\cdots
		   \,+\, x^{-m+3}\raisebox{-2mm}
				{\epsfxsize0.3in\epsffile{basicm-1.ps}}
		   \,+\, x^{-m+1}\raisebox{-2mm}
				{\epsfxsize0.3in\epsffile{basicm.ps}}\\[5pt]
		  &=& x^{m-1} A_{0,m-1}\\
		  &&\,+\,x^{-1} \left(
				x^{m-2}\raisebox{-2mm}
				{\epsfxsize0.3in\epsffile{basic2.ps}}
		   \,+\, x^{m-4}\raisebox{-2mm}
				{\epsfxsize0.3in\epsffile{basic3.ps}}
		   \,+\,\cdots
		   \,+\, x^{-m+4}\raisebox{-2mm}
				{\epsfxsize0.3in\epsffile{basicm-1.ps}}
		   \,+\,x^{-m+2}\raisebox{-2mm}
				{\epsfxsize0.3in\epsffile{basicm.ps}}
		   \,\right)\\[5pt]
		&=&x^{m-1} A_{0,m-1}+(s-s^{-1})P_{m-1} A_{0,0}\\
	       	& &\quad+\,x\left(
				x^{m-2}\raisebox{-2mm}
		 	        {\epsfxsize0.3in\epsffile{basic1.ps}}
		   \,+\, x^{m-4}\raisebox{-2mm}
		 	        {\epsfxsize0.3in\epsffile{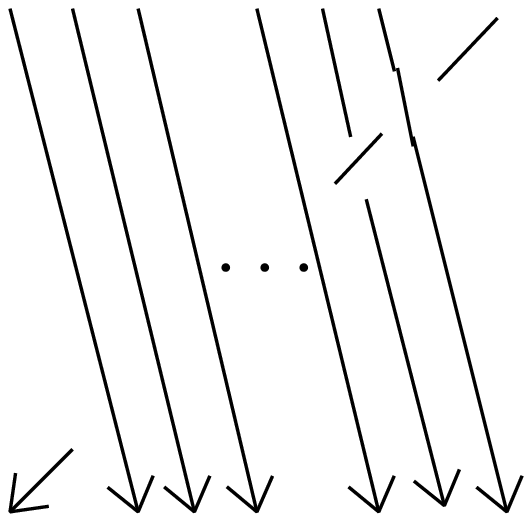}}
		   \,+\,\cdots
		   \,+\, x^{-m+4}\raisebox{-2mm}
		 	        {\epsfxsize0.3in\epsffile{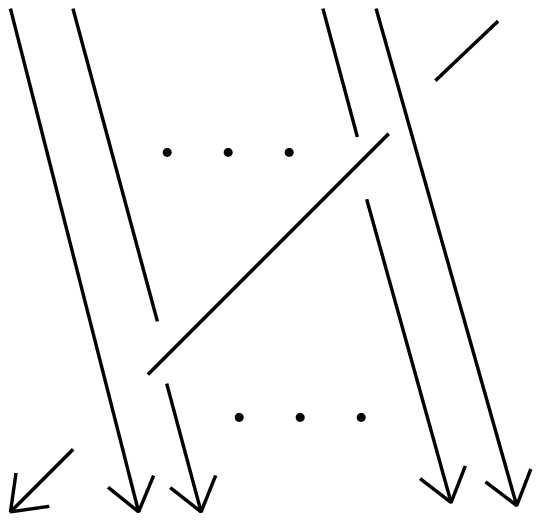}}
		   \,+\, x^{-m+2}\raisebox{-2mm}
		 	        {\epsfxsize0.3in\epsffile{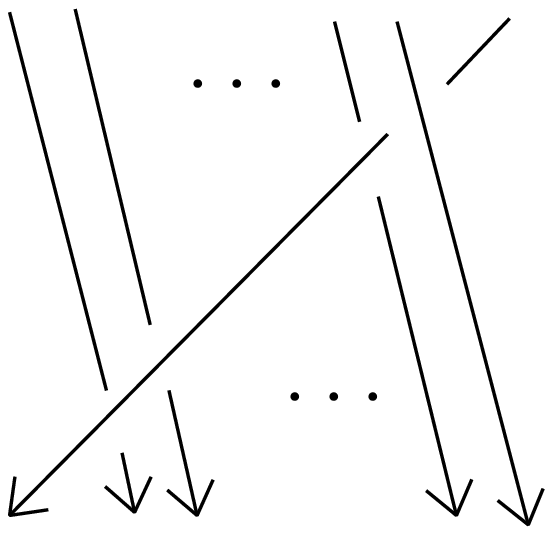}}
		   \,\right)\\[5pt]
		&=&2x^{m-1} A_{0,m-1}+(s-s^{-1})P_{m-1} A_{0,0}\\
		& &\quad+\, x^{m-3}\raisebox{-2mm}
		 	        {\epsfxsize0.3in\epsffile{topneg1.ps}}
		   \,+\,\cdots
		   \,+\, x^{-m+5}\raisebox{-2mm}
		 	        {\epsfxsize0.3in\epsffile{topnegm-2.ps}}
		   \,+\, x^{-m+3}\raisebox{-2mm}
		 	        {\epsfxsize0.3in\epsffile{topnegm-1.ps}}\,\,.\\
\end{eqnarray*}
Applying the skein relation to the $(m-2)$nd crossing
we then see that
\begin{eqnarray*}
P_{m}&=&2x^{m-1} A_{0,m-1}+(s-s^{-1})(P_{m-1} A_{0,0}
	+xP_{m-2} A_{0,1})
         +x^{m-1} A_{0,m-1}\\
	& &+\mbox{\, weighted sum of diagrams with the $(m-1)$st and}\\
	& &\quad \mbox{ $(m-2)$nd crossings negative.}\\
\end{eqnarray*}
Applying the skein relation to all the positive $(m-3)$rd,  
$(m-4)$th, $\ldots$, $2$nd crossings we finally arrive at 
\begin{eqnarray*}
P_m	&=&(m-1)x^{m-1} A_{0,m-1}
			+(s-s^{-1})\sum^{m-1}_{k=2} x^{m-1-i}P_k A_{0,m-k-1} 
			+x^{m-3}\raisebox{-2mm}{\epsfxsize0.3in\epsffile{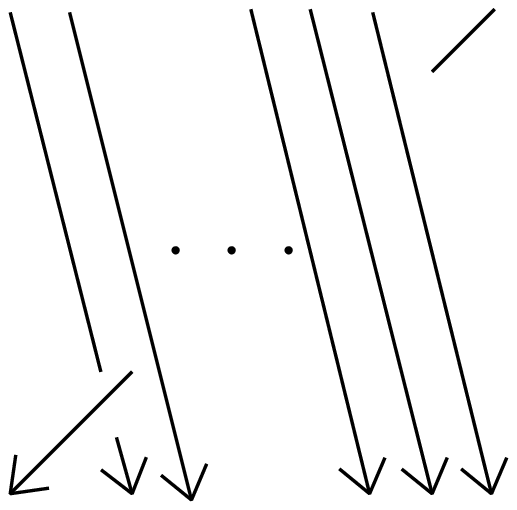}}\\
	&=&(m-1)x^{m-1} A_{0,m-1}
		+(s-s^{-1})\sum^{m-1}_{k=2} x^{m-1-i}P_k A_{0,m-k-1}\\
	&\phantom{a}&\phantom{asdfghjkl;qwertyuiopasdfghkjl;asdfgh}
			+x^{-1}\left(x^{m-2}\raisebox{-2mm}
				{\epsfxsize0.3in\epsffile{basic8.ps}}
			\,\right)\\
	   &=&(m-1)x^{m-1} A_{0,m-1}
		+(s-s^{-1})\sum^{m-1}_{k=2} x^{m-1-k}P_k A_{0,m-k-1}\\
	   & &\quad+(s-s^{-1})P_1 A_{o,m-2}+xx^{m-2} A_{o,m-1}\\[5pt]
	   &=&mx^{m-1} A_{0,m-1}
		+(s-s^{-1})\sum^{m-1}_{k=1} x^{m-1-k}P_k A_{0,m-k-1}\\	
\end{eqnarray*}
\end{proof}

\subsubsection{Theorem.}
\label{xbiff}

The image of the $m$th Adams operator of the fundamental 
representation in the skein of the annulus ${\cal C}^+$, under $\theta$,
is a scalar multiple of
\[
P_m\,=\,x^{m-1}\,\raisebox{-2mm}{\epsfxsize0.3in\epsffile{basic1.ps}}
		   \,+\,x^{m-3}\,\raisebox{-2mm}
				{\epsfxsize0.3in\epsffile{basic2.ps}}
		   \,+\,x^{m-5}\,\raisebox{-2mm}
				{\epsfxsize0.3in\epsffile{basic3.ps}}
		   \,+\,\cdots\,
		     +\, x^{-m+3}\,\raisebox{-2mm}
				{\epsfxsize0.3in\epsffile{basicm-1.ps}}
		   \,+\,x^{-m+1}\raisebox{-2mm}
				{\epsfxsize0.3in\epsffile{basicm.ps}}
\]
namely, 
\[
P_m=[m]\theta(\psi_m( c_1))\,.
\]
\begin{proof}
The proof goes by induction on $m$.
For $m=1$ we have $P_1= c_1$ and $\psi_1( c_1)= c_1$.  

Now, assuming the result holds for all $k<m$, 
we can apply Lemma~\ref{mark} to obtain
the following,
\[
P_m=\sum_{k=1}^{m-1}(s-s^{-1})x^{m-1-k}P_k A_{0,m-1-k}
					\quad+mx^{m-1} A_{0,m-1}\,.
\]
Therefore, by the induction hypothesis,
\[
P_m=\sum_{k=1}^{m-1}
	(s^k-s^{-k})x^{m-1-k}\theta(\psi_k(c_1)) A_{0,m-1-k}
\quad+mx^{m-1} A_{0,m-1}\,.
\]
Using the definition in Lemma~\ref{rich} and Eq.~(\ref{lucy}) 
it can be shown that 
\[
\sum_{k=1}^{m-1}x^{m-1-k}(s^k-s^{-k})\theta(\psi_k(c_1)) A_{0,m-k}
\] 
is the coefficient of $X^{m-2}$ in 
\[
s\,\theta\left(\,\Psi(sX)^{\phantom{1}}\!\right)\Phi^-(xX) 
		- s^{-1}\,\theta\left(\,\Psi(s^{-1}X)\,\right)\Phi^-(xX)\,.
\]
Applying Lemma~\ref{rich} and Prop.~\ref{mary} 
\begin{eqnarray*}
s\theta\left(\,\Psi(sX)\right)\Phi^-(xX)
	&-&s^{-1}\theta\left(\,\Psi(s^{-1}X)\right)\Phi^-(xX)\\
	&=&s\,\theta\left(-C'(sX)D(sX)C(sx^{-1}xX){D'}_q(x^{-1}xX)\right)\\
&\,&\!+s^{-1}\theta\left(D'(s^{-1}X)C(s^{-1}X)
				D(x^{-1}s^{-1}xX){C'}_q(x^{-1}xX)\right)\\
	&=&\!\!\theta\left(s^{-1}D'(s^{-1}X){C'}_q(X)-s{D'}_q(X)C'(sX)\right)
					\,\mbox{by Prop.\ref{equal1}.}
\end{eqnarray*}
We will denote the preimage of the coefficient of $X^{m-2}$ 
in this expression by $g_{m-2}$.  
\[
g_{m-2}=\sum^{m-1}_{k=1}(-1)^k [k](m-k)s^{k-m}c_k d_{m-k}
\ -\sum^{m-1}_{k=1}(-1)^k[m-k]k s^k c_kd_{m-k}.
\]
The coefficient of $c_kd_{m-k}$ in $g_{m-2}$ is 
\begin{eqnarray*}
(-1)^k\left((m-k){(s^{-m+2k}-s^{-m})\over s-s^{-1}}\right.&-&\left.k{(s^m-s^{-m+2k})\over s-s^{-1}}\right)\\
	&=&(-1)^k\left(ms^{-m+k}[k]-k{(s^m-s^{-m})\over s-s^{-1}}\right)\\
	&=&(-1)^k(ms^{-m+k}[k]-k[m])\\
	&=&(-1)^kms^{-m+k}[k]+(-1)^{k-1}k[m]\,.\\ 
\end{eqnarray*}
Recall that $ A_{0,m-1}= \overline{A}_{m}$
and add back on the term $mx^{m-1}\theta( A_{0,m-1})$. 
By Prop.~\ref{karen},
\begin{eqnarray*}
P_m&=&\theta\left(\sum_{k=1}^{m-1}(-1)^k(ms^{-m+k}[k]-k[m])c_kd_{m-k}\right.\\
   &\phantom{a}&\phantom{asdfghjklqwertyuiop}\qquad\qquad
		+m\left.\sum_{k=0}^{m-1}(-1)^ks^k[m-k]c_kd_{m-k}\right)\\[5pt]
   &=&\theta\left(\sum_{k=1}^{m-1}(-1)^k\left(ms^k{(s^{k-m}-s^{-k-m}
		+s^{m-k}-s^{k-m})\over s-s^{-1}}-k[m]\right)c_kd_{m-k}\right.\\
   &\phantom{a}&\left.\phantom{\sum_1^1asdfghjklqwertyuiopzxcvbnmqwe}
					\qquad\qquad+m[m]c_0d_m\right)\\[5pt]
   &=&\theta\left(\sum_{k=1}^{m-1}(-1)^k(m[m]-k[m])c_kd_{m-k}\quad
					+\quad m[m]c_0d_{m}\right)\\[5pt]
   &=&\theta\left([m]\sum_{k=1}^m(-1)^{k-1}kc_kd_{m-k}
				\quad-\quad(-1)^{m-1}m[m]c_md_0\right.\\
   & &\phantom{asdfghkjklqe}\qquad\qquad+m[m]c_0d_m\quad+
		m[m]\left.\sum_{k=1}^{m-1}(-1)^{k}c_kd_{m-k}\right)\\[5pt]
   &=&\theta\left([m]\psi_m(c_1)+m[m]\sum_{k=0}^m(-1)^kc_kd_{m-k}\right)
		\qquad\mbox{by Prop.~\ref{sum}}\\[5pt]
   &=&[m]\theta\left(\psi_m(c_1)\right)+m[m]0\qquad
				\mbox{by Prop.~\ref{equal1}}\\[5pt]
   &=&[m]\theta\left(\psi_m(c_1)\right)\\
\end{eqnarray*}
\end{proof}

\subsection{Vassiliev Invariants.}
\label{sec7}

One motivation for considering the Adams operators came from the work
of Bar-Natan \cite{barnatan} in the context of Vassiliev invariants.
He defines an ``Adams operator'' on chord diagrams.  
We will denote the Adams operator defined
on the subspace of $n$-chord diagrams by $\psi_m^{(n)}$.
Here we investigate the extent
to which applying our Adams operators to the colour in a quantum invariant
upstairs, agrees with Bar-Natan's Adams operators on chord diagrams.
For a detailed account of the theory of Vassiliev invariants the reader is
referred to \cite{barnatan}.

A fundamental result in the theory of Vassiliev invariants
is that every degree $n$ weight system comes from a type $n$ Vassiliev 
invariant and this is well defined up to Vassiliev invariants of lower 
order.  Further, there is a universal Vassiliev invariant $Z$ satisfying
$V_n(W)=W\circ Z$ for each degree $n$ weight system $W$, so we can 
``integrate'' a weight system to get a Vassiliev invariant.  In fact, there 
is a preferred choice for $Z$, which we will denote ${\bf Z}^{{\bf K}}$.  
This is the Kontsevich integral and there are several 
equivalent constructions for ${\bf Z}^{{\bf K}}$, for example \cite{cass, lemur}.
Let ${\bf Z}^{{\bf K}}_n$ be the projection of ${\bf Z}^{{\bf K}}$ onto the 
$n$th graded
piece of the algebra of chord diagrams and suppose that $V$ is a type $n$
Vassiliev invariant.  We say that $V$ is \defn{canonical} if and only if 
$V=W(V)\circ{\bf Z}^{{\bf K}}_n$.  
Let $h$ be a formal parameter and suppose $V=\sum_{i=0}^\infty V_i h^i$
is a power series over $h$ where $V_i$ is a Vassiliev invariant of type $i$.
We say $V$ is \defn{canonical} if each of the $V_i$ is a canonical invariant.

We are interested in canonical invariants for the following reasons.  
Two canonical invariants are equal if and only if
their weight systems are equal.  Equivalently, the weight system uniquely
determines the canonical invariant.
The $U_q(sl(N))$-invariant $J(K;V_\lambda)$ with $\lambda$ fixed is a canonical 
power series, when expanded at $q=e^h$ \cite{cass, lemur}.

Let $D$ be an $n$-chord diagram. Following \cite[Defn. 3.11]{barnatan}, we 
define $\psi_m^{(n)}(D)$ to be the
sum of all possible ways of lifting $D$ to the $m$th cover of the circle. 
For example,
\begin{eqnarray*}
\psi_2^{(2)}\,\left(\ \raisebox{-3mm}{\epsfxsize.3in\epsffile{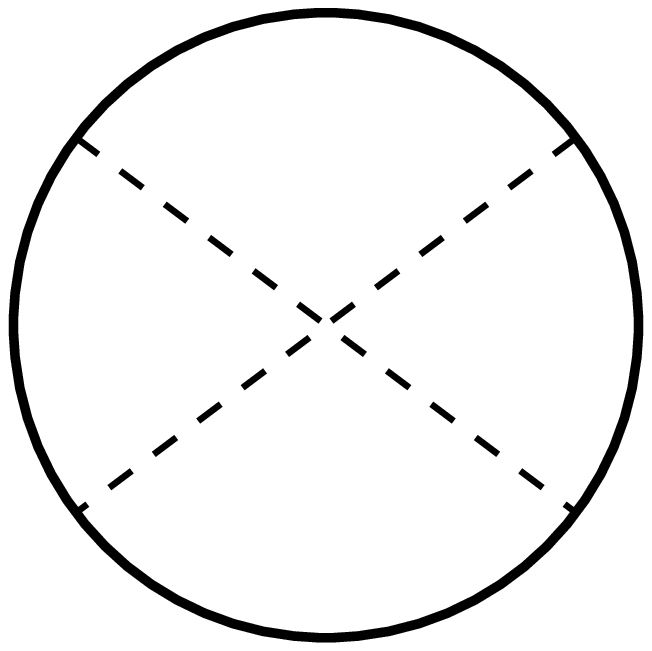}}
		\right)
\quad&=&\quad\raisebox{-4mm}{\epsfxsize.4in\epsffile{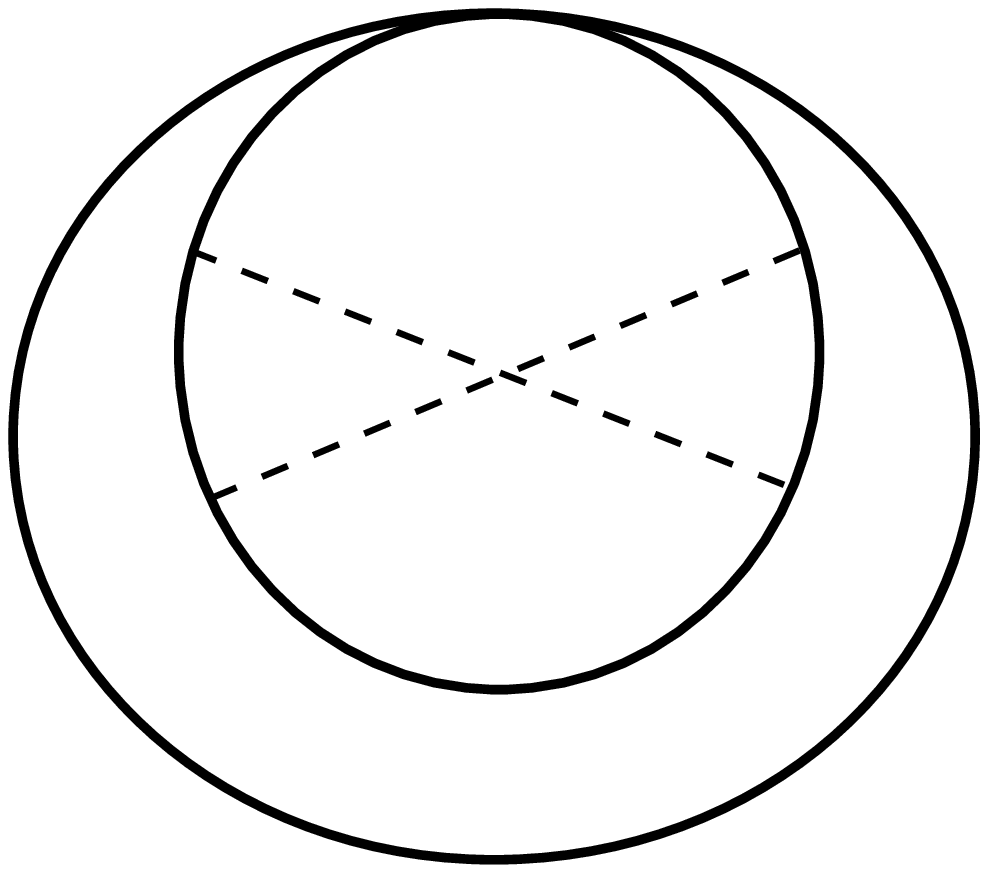}}
\quad+\quad\raisebox{-4mm}{\epsfxsize.4in\epsffile{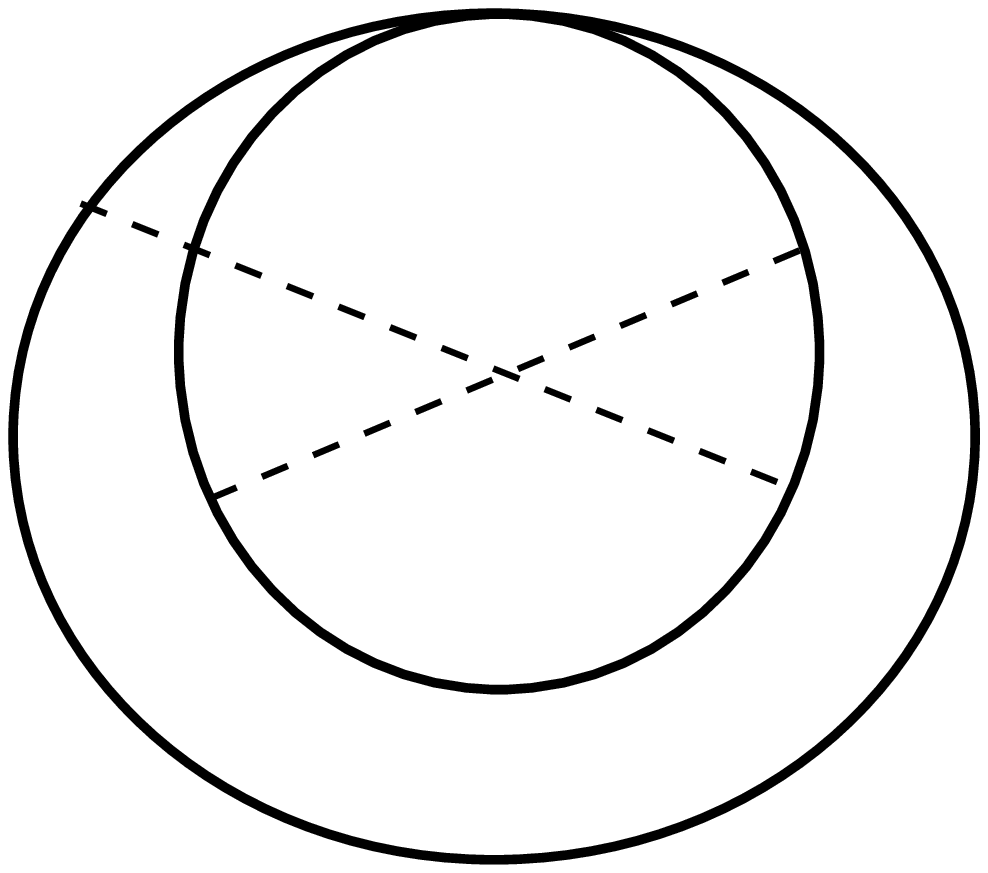}}
\quad+\quad\raisebox{-4mm}{\epsfxsize.4in\epsffile{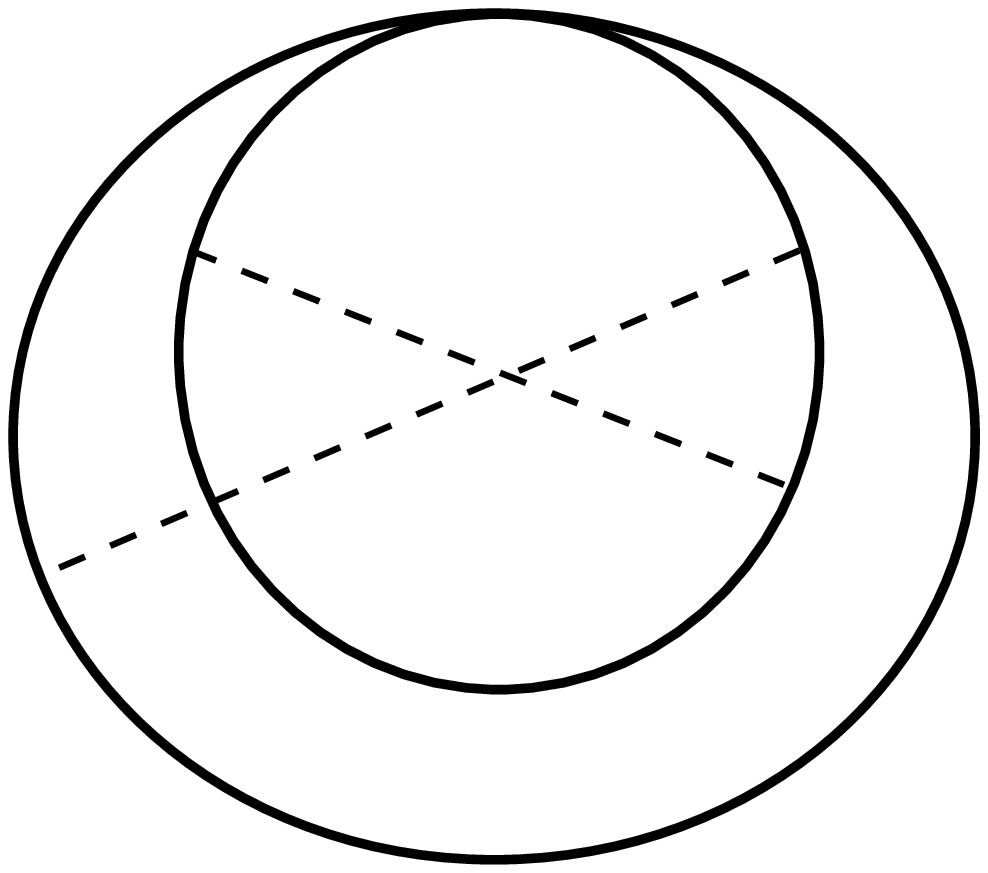}}\quad + \quad \cdots\\
&=&8\ \raisebox{-3mm}{\epsfxsize.3in\epsffile{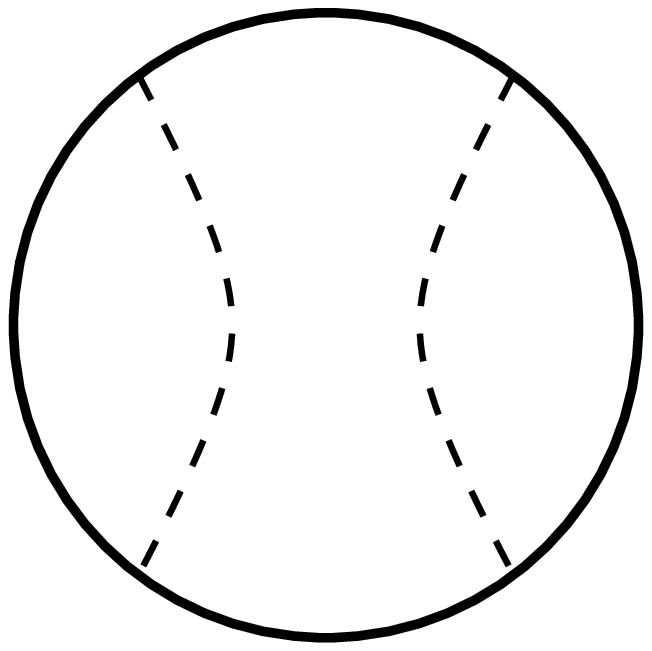}}
\quad+\quad 8\ \raisebox{-3mm}{\epsfxsize.3in\epsffile{chord.ps}}\,.
\end{eqnarray*}

\subsubsection{Theorem.} 

Let $J(K;c_1)=\sum_{i=0}^\infty V_i h^i$
be the Vassiliev power series expansion of 
$J(K; c_1)$. 
Take the Vassiliev power series expansion of $J(K;\psi_m(c_1))$
to be $\sum_{i=0}^\infty U_i(K) h^i$.
Then $U_i$ is the canonical degree $i$ Vassiliev invariant with weight system 
$W_i\circ \psi_m^{(i)}$.
\[
U_i=\left(W(V_i)\circ\psi_m^{(i)}\right)\circ{\bf Z}^{\bf K}_i
\]
\begin{proof}
By Theorem \ref{xbiff}, 
\[
J(K;\psi_m(c_1))\,=\,{1\over [m] }\,J(K*P_m; c_1)\,=\,{1\over [m] }\,
					\sum_{i=0}^\infty V_i(K*P_m) h^i
\]

Assume that $K$ has $n$ double points.
Then $K*P_m$ can be written as a linear combination of singular knots 
with $n$ double points and $V_i(K*P_m;c_1)=0$ when $i<n$. 
Further, expanding $[m]$ as a power series in $h$, $[m]=m+O(h)$.
Thus
\begin{equation}
J(K;\psi_m(c_1))\,=\,{1\over [m]}\,J(K*P_m;c_1)\,=\,{1\over m}\,V_n(K*P_m)h^n+ O(h^{n+1})
\label{haircut}
\end{equation}

Using the notation of Sect.~\ref{sec6}, consider a summand 
$A_{i,m-1-i}$ of $P_m$.  Since $V_n$ is a Vassiliev invariant,
\[
V_n(K*A_{i,m-1-i})= V_n(K*A_{i+1,m-2-i})\,-\,V_n(K*B_i)\,,
\]
where $B_i$ is a pattern which contains a double point. Since 
we assumed that
$K$ has $n$ double points, $V_n(K*B_i)=0$.  By repeated application
of this relation, we see that $V_n(K*A_{i,m-1-i})=V_n(K*A_m)$, where 
$A_m$($=A_{m-1,0}$) is the element depicted in Fig.~\ref{am}. 
Since $x=1+O(h)$, as a power series
in $h$, 
\[
V_n(K*P_m)\,=\,m V_n(K*A_m) +O(h)\,.
\]
Substituting back into Eq.~(\ref{haircut}), we have
\begin{equation}
J(K;\psi_m(c_1))={m\over m} V_n(K*A_m)h^n +O(h^{n+1})\,.
\label{swim}
\end{equation}
Equating coefficients of $h^n$, Eq.~\ref{swim} implies that
for a knot with $n$ double points
$V_n(K*A_m)=U_n(K)$.  Equivalently, for any $n\geq 0$ and any knot $K$,
the two invariants agree at the weight system level,
\[
W(V_n(K*A_m)) = W(U_n(K))\,.
\]
The satellite $K*A_m$ is a connected $m$ string cable of $K$. 
We can think of $V_n(K*A_m)$ as an invariant of $K$
by taking the composition of $V_n$ with
the cabling map $\varphi_m$ : $\varphi_m(K)=K*A_m$.
Bar-Natan \cite{barnatan} showed that 
$W(V_n\circ \varphi_m)=W(V_n)\circ\psi_m^{(n)}$
and as we commented above, $U_n$ is a canonical invariant,
therefore, for each $n\geq 0$,
\[
U_n=\left(\,W(V_n)\circ\psi_m^{(n)}\,\right)\circ {\bf Z}^{{\bf K}}_n\,.
\]
\end{proof}

We have shown that the canonical Vassiliev power series 
with weight system $W(J(K;c_1))\circ \psi_m^{(n)}$ is the expansion
of $J(K;\psi_m(c_1))$. 
(By $W(J(K;c_1))$ we mean the sum of the weight systems for each of
the Vassiliev invariants in the Vassiliev power series expansion
of $J(K;c_1)$.)
Bar-Natan's result holds in a more general context.  It says that 
for any type $n$ Vassiliev invariant $V$, 
\[
W(V\circ\varphi_m)=W(V)\circ\psi_m^{(n)}\,.
\]
In particular, if we can express $\psi_m(\lambda)$ as a linear combination
of cabling patterns, coloured by $\lambda$, we could prove the following 
conjecture

\subsubsection{Conjecture.}

Let $J(K;V_\lambda)$ have a Vassiliev power series expansion 
$\sum_{i=0}^\infty V_{i,\lambda} h^i$.  Then the 
canonical Vassiliev invariant with weight system 
$W(V_{i,\lambda})\circ\psi_m^{(n)}$ is the coefficient of $h^i$ in
the Vassiliev power series expansion of $J(K;\psi_m(\lambda))$.
\qed

An obvious guess is that if we colour $P_m$ (or some other linear 
combination of the braids in $P_m$) with $\lambda$ this might
work.  However, this doesn't work even for the simplest case, 
$\psi_2(c_2)$. In terms of Young diagrams, 
\[
\psi_2(c_2)=(4)-(2,1,1)+(2,2)\,.
\]
We would need the following equation to hold in ${\cal C}^+$,
\[
Q_{(4)} - Q_{(2,1,1)} + Q_{(2,2)} 
\ =\ 
A\,\raisebox{-.5cm}{\epsfxsize.6in\epsffile{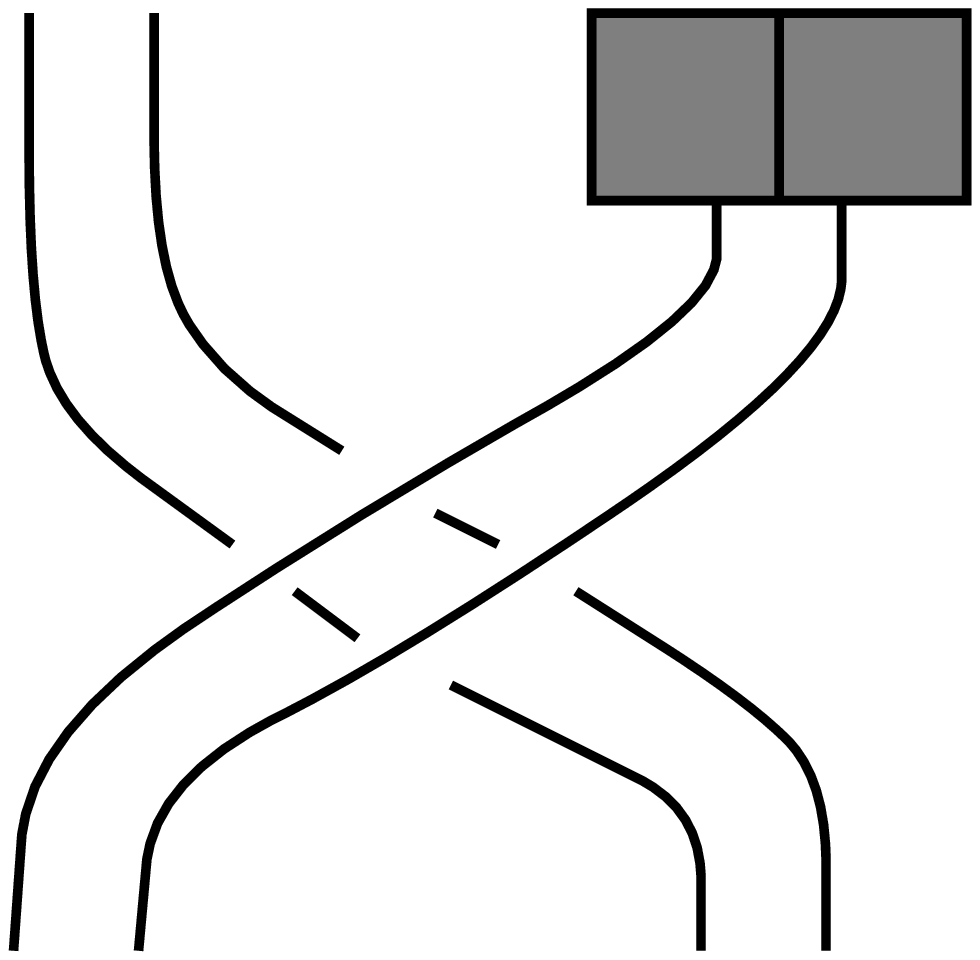}}\ 
+\ B\,\raisebox{-.5cm}{\epsfxsize.6in\epsffile{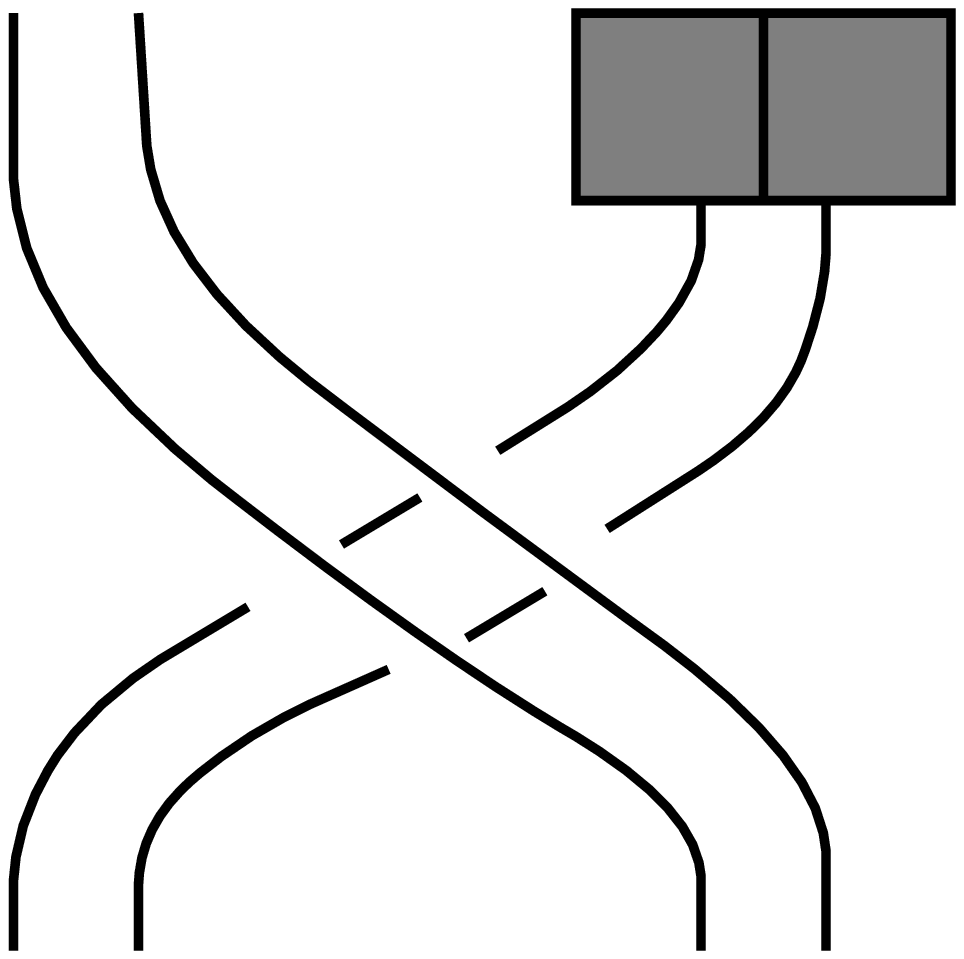}}\,.
\]
Expressing everything in terms of monomials in $\{A_m\}_{m\in \NS}$, 
we can compare the coefficients of each of the five weighted degree four
monomials in the $A_m$, and so solve 
for $A$ and $B$ in more than one way.
Unfortunately, the answers are not consistent.   
Thus, there is no linear combination of connected $2$-cables, $P'$, for which
$J(K;\psi_2(c_2))=J(K*P';c_2)$.  
The pattern for $\psi_2(c_2)$ must, therefore, involve the $2$-parallel
and make use of the Vassiliev skein relation, to relate it to the connected 
cablings  modulo extra double points.  Since $(2,2)$ is not a hook shaped
diagram, perhaps this is not too surprising in light of Cor.~\ref{sparks}.

\bibliography{library}
\end{document}